\documentclass[aap,MSNbibl,citesort,dvips]{arximspdf}
\usepackage{mathbh}
\usepackage{mathrsfs}

\doi{10.1214/12-AAP843} 
\volume{23}
\issue{5}
\pubyear{2013}
\firstpage{1721}
\lastpage{1754}

\makeatletter
\newcommand{\eqref}[1]{(\ref{#1})}
\newcommand{\Dc}{\mathcal{D}}

\newcommand{\indic}{\mathbh{1}}
\newcommand{\X}{\mathbb{X}}
\def\bm{\mathbf}
\def\mathds{\mathbb}
\newcommand{\E}{\mathds{E}}
\newcommand{\N}{\mathds{N}}
\renewcommand{\P}{\mathds{P}}

\newcommand{\Pcr}{\mathscr{P}}
\newcommand{\Xcr}{\mathscr{X}}
\newcommand{\Ccr}{\mathscr{C}}

\newcommand{\ddr}{\mathrm{d}}

\newtheorem{thmm}{Theorem}
\newproclaim{defi}{Definition}
\newtheorem{lem}{Lemma}
\newtheorem{cor}{Corollary}
\newproclaim{rmk}{Remark}
\newtheorem{prp}{Proposition}
\newproclaim{exe}{Example}
\makeatother

\begin{document}
\begin{frontmatter}

\title{Conditional formulae for Gibbs-type exchangeable random
partitions\thanksref{T1}}
\runtitle{Conditional formulae for Gibbs-type partitions}
\thankstext{T1}{Supported by the European Research Council (ERC) through StG ``N-BNP'' 306406.}

\begin{aug}
\author[A]{\fnms{Stefano}~\snm{Favaro}\ead[label=e1]{stefano.favaro@unito.it}\thanksref{t2}},
\author[B]{\fnms{Antonio}~\snm{Lijoi}\corref{}\ead[label=e2]{lijoi@unipv.it}\thanksref{t2}}
\and
\author[A]{\fnms{Igor}~\snm{Pr\"unster}\ead[label=e3]{igor@econ.unito.it}\thanksref{t2}}
\thankstext{t2}{Also affiliated to Collegio Carlo Alberto, Moncalieri, Italy.}

\runauthor{S. Favaro, A. Lijoi and I. Pr\"unster}

\affiliation{University of Torino, University of Pavia and University
of Torino}
\address[A]{S. Favaro\\
I. Pr\"unster\\
Department of Economics and Statistics\\
University of Torino\\
Corso Unione Sovietica 218/bis, I-10134 Torino\\
Italy\\
\printead{e1}\\
\phantom{E-mail: }\printead*{e3}}

\address[B]{A. Lijoi\\
Department of Economics and Business\\
University of Pavia\\
Via San Felice 5, 27100 Pavia\\
Italy \\
\printead{e2}}

\end{aug}

\received{\smonth{6} \syear{2011}}
\revised{\smonth{1} \syear{2012}}

%
\begin{abstract}
Gibbs-type random probability measures and the exchangeable random
partitions they induce represent an important framework both from a
theoretical and applied point of view. In the present paper, motivated
by species sampling problems, we investigate some properties concerning
the conditional distribution of the number of blocks with a certain
frequency generated by Gibbs-type random partitions. The general
results are then specialized to three noteworthy examples yielding
completely explicit expressions of their distributions, moments and
asymptotic behaviors. Such expressions can be interpreted as Bayesian
nonparametric estimators of the rare species variety and their
performance is tested on some real genomic data.
\end{abstract}

%
\begin{keyword}[class=AMS]
\kwd{60G57}
\kwd{62G05}
\kwd{62F15}.
\end{keyword}

\begin{keyword}
\kwd{Bayesian nonparametrics}
\kwd{Exchangeable random partitions}
\kwd{Gibbs-type random partitions}
\kwd{sampling formulae}
\kwd{small blocks}
\kwd{species sampling problems}
\kwd{$\sigma$-diversity}.
\end{keyword}

\end{frontmatter}

\section{Introduction}\label{sec1} Let $\mathbb{X}$ be a complete and separable
metric space equipped with the Borel $\sigma$-algebra $\mathscr X$ and
denote by $\mathscr P$ the space of probability distributions defined
on $(\X,\Xcr)$ with $\sigma(\mathscr P)$ denoting the Borel $\sigma
$-algebra of subsets of $\Pcr$. By virtue of de~Finetti's
representation theorem, a sequence of $\X$-valued random elements
$(X_{n})_{n\geq1}$, defined on some probability space $(\Omega
,\mathscr
{F},\mathds{P})$, is exchangeable if and only if there exists a
probability measure $Q$ on the space of probability distributions
$(\mathscr P, \sigma(\mathscr P))$ such that
%
%
\begin{equation}\label{eqdefin}
\P[X_1\in A_1,\ldots,X_n\in A_n]=\int_{\mathscr P}\prod_{i=1}^n
P(A_i)  Q(\ddr P)
\end{equation}
for any $A_1,\ldots,A_n$ in $\Xcr$ and $n\ge1$. The probability
measure $Q$ directing the exchangeable sequence $(X_{n})_{n\geq1}$ is
also termed \textit{de Finetti measure} and takes on the interpretation
of prior distribution in Bayesian applications. The representation
theorem can be equivalently stated by saying that, given an
exchangeable sequence $(X_{n})_{n\geq1}$, there exists a random
probability measure (r.p.m.) $\tilde{P}$, defined on $(\X,\Xcr)$ and
taking values in $(\mathscr P, \sigma(\mathscr P))$, such that
%
%
\begin{equation}
\label{eqexchangeabl}
\P[X_{1}\in A_{1},\ldots,X_{n}\in A_{n}|\tilde P]=\prod
_{i=1}^{n}\tilde
{P}(A_{i})
\end{equation}
almost surely, for any $A_{1},\ldots,A_{n}$ in $\mathscr X$ and $n\geq
1$. In this paper we will focus attention on almost surely discrete
r.p.m.'s, that is, $\tilde P$ is such that $\P[\tilde P\in\mathscr
P_d]=1$
with $\mathscr P_d\in \sigma(\mathscr P)$ indicating the set
of discrete probability
measures on $(\X,\Xcr)$ or, equivalently, $(X_n)_{n\ge1}$ is directed
by a de Finetti measure $Q$ that is concentrated on~$\mathscr P_d$. An
almost surely discrete r.p.m.~(without fixed atoms) can always be
written as
%
%
\begin{equation}
\tilde P=\sum_{i \ge1} \tilde p_i \delta_{\hat X_i}
\label{eqdiscrete}
\end{equation}
for some sequences $(\hat{X}_i)_{i\ge1}$ and $(\tilde p_i)_{i\ge1}$
of, respectively, $\X$-valued random locations and nonnegative
random weights such that $\P[\sum_{i\ge1}\tilde p_i=1]=1$ almost
surely.

In the following we will assume that the two sequences in \eqref
{eqdiscrete} are independent. These specifications imply that a sample
$(X_1,\ldots,X_n)$ from the exchangeable sequence generates a random
partition $\Pi_n$ of the set of integers $\N_n:=\{1,\ldots,n\}$, in the
sense that any $i\ne j$ belongs to the same partition set if and only
if $X_i=X_j$. The random number of partition sets in $\Pi_n$ is denoted
as $K_n$ with respective frequencies $N_{1},\ldots,N_{K_n}$.
Accordingly, the sequence $(X_n)_{n\ge1}$ associated to a
r.p.m. $\tilde P$ as in \eqref{eqdiscrete} induces an exchangeable
random partition $\Pi=(\Pi_n)_{n\ge1}$ of the set of natural numbers
$\N$. The distribution of $\Pi$ is characterized by the sequence $\{p_k^{(n)}\dvtx   1\le k\le n,  n\ge1\}$ such that
%
%
\begin{equation}
\label{eqeppfdef}
p_k^{(n)}(\bm{n})=\P[\Pi_n=\pi],
\end{equation}
for $\pi$ a partition of $\mathbb{N}_n$ into $k$ blocks with vector frequencies $\mathbf{n} = (n_1,\ldots, n_k)$ such that
$\sum^k_{j=1} n_j = n$.
Hence, \eqref{eqeppfdef} identifies, for any $n\ge1$, the
probability distribution of the random partition $\Pi_n$ of $\N_n$ and
is known as the \textit{exchangeable partition probability function}
(EPPF), a concept introduced by J.~Pitman \cite{Pit95} as a major
development of earlier results on exchangeable random partitions due to
J.~F.~C.~Kingman (see, e.g., \cite{Kin78,Kin82}). It is worth noting
that EPPFs can be defined either by starting from an exchangeable
sequence associated to a discrete r.p.m.~and looking at the induced
partitions or by defining directly the partition distribution. In the
latter case, the distribution of the random partitions $\Pi_n$ must
satisfy certain consistency conditions and a symmetry property that
guarantees exchangeability. A comprehensive account on exchangeable
random partitions can be found in \cite{Pit06} together with an
overview of the numerous application areas and relevant references.

\subsection{Gibbs-type r.p.m.'s and partitions}\label{sec1.1} We now recall the
definition of a general class of r.p.m.'s and of the exchangeable random
partitions they induce together with some of distinguished special
cases. This important class, introduced and thoroughly studied in \cite
{Gne05}, is characterized by the fact that its members induce
exchangeable random partitions admitting EPPFs with product form, a
feature which is crucial for guaranteeing mathematical tractability.
Before introducing the definition, set $\mathcal{D}_{n,j}:=\{
(n_{1},\ldots,n_{j})\in\{1,\ldots,n\}^{j}\dvtx   \sum_{i=1}^{j}n_{i}=n\}$
and denote by $(a)_{q}=\Gamma(a+q)/\Gamma(a)$ the $q$th ascending
factorial of $a$ for any integer $q \geq 1$.

\begin{defi}\label{gibbspartition}
Let $(X_n)_{n\geq1}$ be an exchangeable sequence associated to an
almost surely discrete r.p.m.~\eqref{eqdiscrete} for which locations
$(\hat{X}_i)_{i\ge1}$ and weights $(\tilde p_i)_{i\ge1}$ are
independent. Then the r.p.m.~$\tilde P$ and the induced exchangeable
random partition are said of \textit{Gibbs-type} if, for any $n\ge1$,
$1\leq j\leq n$ and $(n_{1},\ldots,n_{j})\in\Dc_{n,j}$, the
corresponding EPPF can be represented as follows:
%
%
\begin{equation}\label{eqgibbseppf}
p_{j}^{(n)}(n_{1},\ldots,n_{j})=V_{n,j}\prod
_{i=1}^{j}(1-\sigma)_{n_{i}-1}
\end{equation}
for $\sigma\in(-\infty, 1)$ and a set of nonnegative weights $
\{
V_{n,j}\dvtx n\geq1, 1\leq j\leq n\}$ satisfying the recursion
$V_{n,j}=V_{n+1,j+1}+(n-\sigma j)V_{n+1,j}$ with $V_{1,1}=1$.
\end{defi}

Hence, a Gibbs-type random partition is completely specified by the
choice of $\sigma<1$ and the weights $V_{n,j}$'s. The role of $\sigma$
is crucial since it determines the clustering structure as well as the
asymptotic behavior of Gibbs-type models. As for the latter aspect,
for any $n\ge1$ define
\[
c_{n}(\sigma):=\indic_{(-\infty,0)}(\sigma)+\log(n) \indic_{\{0\}
}(\sigma)+n^\sigma \indic_{(0,1)}(\sigma)
\]
where $\indic_A$ denotes the indicator function of set $A$.
Then, for any Gibbs-type r.p.m.~there exists a strictly positive and
almost surely finite random variable $S_{\sigma}$, usually termed
$\sigma$-\textit{diversity}, such that
%
%
\begin{equation}\label{eqsdiversity}
\frac{K_{n}}{c_{n}(\sigma)}\stackrel{\mathrm{a.s.}}{\longrightarrow
}S_{\sigma},
\end{equation}
for $n\rightarrow+\infty$. See \cite{Pit03}, Section 6.1, for details.
Finally, it is worth recalling that the solutions of the backward
recursions defining the $V_{n,j}$'s form a convex set whose extreme
points are determined in \cite{Gne05}, Theorem~12, providing a
complete characterization of Gibbs-type models according to the values
of $\sigma$ they assume. In the next subsection we concisely point out
three important explicit special cases to be dealt with also in the sequel.


\subsection{Examples}\label{sec1.2} We will illustrate three noteworthy examples of
Gibbs-type r.p.m.'s that correspond to different choices of $\sigma$
and the $V_{n,j}$'s in Definition~\ref{gibbspartition}. The first one
is the \textit{Dirichlet process} \cite{Fer73}, which corresponds to a
Gibbs-type r.p.m.~characterized by $\sigma=0$ and $V_{n,j}=\theta
^j/(\theta)_n$ with $\theta>0$. The implied EPPF coincides with
%
%
\begin{equation}
p_j^{(n)}(n_1,\ldots,n_j)=\frac{\theta^j}{(\theta)_n}\prod_{i=1}^j (n_i-1)!
\label{eqeppfewens}
\end{equation}
and is well known in Population Genetics as the \textit{Ewens model}. See
\cite{Ewe72} and
references therein.

The most interesting special case for our purposes is a generalization
of \eqref{eqeppfewens} that has been provided by J.~Pitman in \cite
{Pit95}. It corresponds to the exchangeable random partition
generated by the \textit{two-parameter Poisson--Dirichlet process},
which coincides with a Gibbs-type r.p.m.~with $\sigma\in(0,1)$ and,
for any $\theta>-\sigma$, $V_{n,j}=\prod_{i=0}^{j-1}(\theta+i\sigma
)/(\theta)_{n}$. The EPPF turns out to be
%
%
\begin{equation}
\label{eqeppf}
p_j^{(n)}(n_1,\ldots,n_j)
=\frac{\prod_{i=0}^{j-1}(\theta+i\sigma) }{(\theta)_{n}}\prod
_{i=1}^{j}(1-\sigma)_{n_{i}-1}.
\end{equation}
Clearly, the Ewens model \eqref{eqeppfewens} is recovered from \eqref
{eqeppf} by letting $\sigma\to0$. The r.p.m. and the partition
distribution associated to \eqref{eqeppf} will be equivalently termed
the $\operatorname{PD}(\sigma,\theta)$ process or \textit{Pitman model}.

Finally, another notable example of Gibbs-type r.p.m.~has been
recently provided in \cite{Gne10}. It is characterized by $\sigma=-1$
and weights of the form
%
%
\begin{equation}\label{eqgnedinweight}
V_{n,j}=(\gamma)_{n-j}\frac{\prod_{i=1}^{j-1}(i^{2}-\gamma i+\zeta
)}{\prod_{i=1}^{n-1}(i^{2}+\gamma i+\zeta)},
\end{equation}
where $\zeta$ and $\gamma$ are chosen such that $\gamma\geq0$ and
$i^{2}-\gamma i+\zeta>0$ for all $i\geq1$. In the sequel we will term
both the r.p.m.~and the induced exchangeable random partition as the
\textit{Gnedin model}.

\subsection{Aims and outline of the paper}\label{sec1.3}
The main applied motivation of the present study is related to species
sampling problems. Indeed, in many applications that arise, for
example, in population genetics, ecology and genomics, a population is
a~composition of individuals (e.g., animals, plants or genes) of
different species: the $\hat{X}_i$'s and the $\tilde p_i$'s in \eqref
{eqdiscrete} can then be seen as species labels and species
proportions, respectively. In most cases one is interested in the
$\tilde p_i$'s or in some functionals of them: this naturally leads to
work with the random partitions induced by an exchangeable sequence.
The number of distinct partition blocks $K_n$ takes on the
interpretation of the number of different species detected in the
observed sample $(X_1,\ldots,X_n)$ and the $N_{j}$'s are the species
frequencies. Given the relevance and intuitiveness of such an applied
framework, throughout the paper we will often resort to the species
metaphor even if the tools we will introduce and the results we will
achieve are of interest beyond the species sampling framework.

Our first goal consists in analyzing certain distributional properties
of Gibbs-type r.p.m.'s. Specifically, we are interested in determining
the probability distribution of the number of partition blocks having a
certain size or frequency. In other words, given an exchangeable
sequence $(X_n)_{n\ge1}$ as in \eqref{eqdefin} associated to a
Gibbs-type r.p.m., we investigate distributional properties of: (i)
the number of species with frequency $l$ in a sample of size $n$,
namely, $M_{l,n}=\sum_{i=1}^{K_n}\indic_{\{l\}}(N_{i})$; (ii) the
number of species $M_{l,n+m}=\sum_{i=1}^{K_{n+m}}\indic_{\{l\}}(N_{i})$
with frequency $l$ in an enlarged sample of size $n+m$, for $m \geq0$,
conditionally on the species composition detected within a $n$-size
sample $(X_1,\ldots,X_n)$. Note that the latter problem is considerably
more challenging since it requires to account for the allocation of
$(X_{n+1},\ldots,X_{n+m})$ between ``old'' and ``new'' species together
with the sequential modification of their frequencies, conditional on
$(X_1,\ldots,X_n)$.

Solving problem (ii) is also the key for achieving our second goal,
namely, the derivation of estimators for \textit{rare species} variety,
where rare species are identified as those with a frequency not greater
than a specific abundance threshold $\tau$. This is of great importance
in numerous applied settings. For example, in ecology conservation of
biodiversity is a fundamental theme and it can be formalized in terms
of the number of species whose frequency is greater than a specified
threshold. Indeed, any form of management on a sustained basis requires
a certain number of sufficiently abundant species (the so-called
breeding stock). We shall address the issue be relying on a Bayesian
nonparametric approach: the de Finetti measure associated to a
Gibbs-type r.p.m.~represents the nonparametric prior distribution and
relying on the conditional (or posterior) distributions in (ii), one
derives the desired estimators as conditional (or posterior) expected
values. Bayesian estimators for overall species variety, namely, the
estimation of the distinct species (regardless of the respective
frequencies), have been introduced and discussed in \cite
{Lij207,lmp07bmc,Lij308,Fav09}. Further contributions at the
interface between Bayesian Nonparametrics and Gibbs-type random
partitions can be found in \cite{griffiths07,Ho,lmp07}. None of the
existing work provides estimators for the number of species with
specific abundance. Here we fill in this important gap and, besides
providing general results valid for the whole family of Gibbs-type
r.p.m.'s, we specialize them to the three examples outlined in
Section~\ref{sec1.2}. This leads to explicit expressions that are of
immediate use in applications.

The paper is structured as follows. Section~\ref{sec2} provides distributional
results on the unconditional structure of $M_{l,n}$ and the conditional
structure of $M_{l,n+m}$, given the species composition detected in a
sample of size $n$, for general Gibbs-type r.p.m.'s together with the
corresponding estimators. Section~\ref{sec3} focuses on the three special cases
of the Dirichlet process, and the models of Gnedin and of Pitman. In
particular, for these special cases we also provide asymptotic results
concerning the conditional distribution of $M_{l,n+m}$, given the
species composition detected in a sample of size $n$, as the size of
the additional sample $m$ increases. The framework for genomic
applications, including platforms under which such estimation problems
arise, is presented in Section~\ref{sec4}, where the methodology is also tested
on real genomic data. In Section~\ref{sec5} the proofs of the results of
Sections~\ref{sec2} and~\ref{sec3} and some useful techniques are described.

\section{Distribution of cluster frequencies}\label{sec2}

\subsection{Probability distribution of $M_{l,n}$}\label{sec2.1}
We start our analysis of distributional properties of Gibbs-type
random partitions by focusing on the unconditional distribution of the
number of blocks with a certain size $l$, $M_{l,n}$. The blocks with
relatively low frequency are typically referred to as \textit{small
blocks} (see, e.g., \cite{Sch10}), which, in terms of species
sampling, will represent the rare species.

First note that the EPPF \eqref
{eqgibbseppf} yields the probability distribution of
$\bm{M}_{n}:=(M_{1,n},\ldots,M_{n,n})$. Specifically, the so-called
Gibbs-type sampling formula determines the probability distribution of
$\bm{M}_{n}$ and it corresponds to
%
%
\begin{equation}\label{eqgibbssampling}
\P[\bm{M}_{n}=(m_{1},\ldots,m_{n})]=V_{n,j}
n! \prod_{i=1}^{n}\biggl(\frac{(1-\sigma)_{i-1}}{i!}
\biggr)^{m_{i}}\frac{1}{m_{i}!},
\end{equation}
for any $(m_1,\ldots,m_n)\in\{0,1,\ldots\}^n$ such that $\sum
_{i=1}^n i
m_i=n$ and $\sum_{i=1}^n m_i=j$. The next proposition provides explicit
expressions for the $r$th factorial moments of $M_{l,n}$ in terms of
generalized factorial coefficients $\Ccr(n,k;\sigma)$. Recall that, for
any $n\geq1$ and $k\leq n$, $\mathscr{C}(n,k;\sigma)$ is defined as
$(\sigma t)_{n}=\break\sum_{k=0}^{n}\mathscr{C}(n,k;\sigma)(t)_{k}$ for
$\sigma\in\mathbb{R}$ and, moreover, is computable as $\mathscr
{C}(n,k;\sigma)=(1/k!)\sum_{j=0}^{k}(-1)^{j}{k\choose j}(-\sigma j)_{n}$
with the proviso $\mathscr{C}(0,0;\sigma)=1$, $\mathscr
{C}(n,0;\sigma
)=0$ for any $n>0$ and
$\mathscr{C}(n,k;\sigma)=0$ for any $k>n$. For an exhaustive account on
generalized factorial coefficients the reader is referred to \cite{Cha05}.

\begin{prp}\label{priormarginal}
Let $(X_n)_{n\ge1}$ be an exchangeable sequence associated to a
Gibbs-type r.p.m. Then, for any $l=1,\ldots,n$ and $r\geq1$,
\begin{equation}\label{eqmomentmarginalprior}
\E\bigl[(M_{l,n})_{[r]}\bigr]=\biggl(\frac{(1-\sigma
)_{l-1}}{l!}
\biggr)^{r}
(n)_{[lr]} \sum_{j=1}^{n}V_{n,j} \frac{\mathscr{C}(n-rl,j-r;\sigma
)}{\sigma^{j-r}},
\end{equation}
where $(a)_{[q]}=a(a-1) \cdots (a-q+1)$ for any $q\ge1$.
\end{prp}

By using standard arguments involving probability generating functions,
one can use the factorial moments \eqref{eqmomentmarginalprior} for
determining the probability distribution of $M_{l,n}$. This will be\vadjust{\goodbreak}
illustrated for the three examples in Section~\ref{sec3}. The asymptotic
behavior of $M_{l,n}$, as $n\to\infty$, is determined in \cite{Pit06}, Lemma~3.11: if $\tilde P$ is a Gibbs-type r.p.m.~with
$\sigma\in(0,1)$, then for any $l\geq1$
%
%
\begin{equation}\label{eqgenerallimit}
\frac{M_{l,n}}{n^{\sigma}} \stackrel{\mathrm{d}}{\longrightarrow}
\frac{\sigma(1-\sigma)_{l-1}}{l!}S_{\sigma}
\end{equation}
as $n \rightarrow+\infty$, where $S_{\sigma}$ is the $\sigma
$-diversity defined in \eqref{eqsdiversity}. Some recent interesting
developments on the asymptotic behavior of the random variable
$M_{l,n}$ associated to a generic exchangeable random partition are
provided in \cite{Sch10}.

\subsection{Conditional formulae}\label{sec2.2} Unlike the study of unconditional
properties of Gibbs-type random partitions, that are the focus of a
well-established literature with plenty of results, the investigation
of conditional properties for this family of partitions has been only
recently started in \cite{Lij308} and many issues are still to be
addressed. We are going to focus on determining the distribution of
$M_{l,n+m}$ conditional on the number of distinct species $K_n$, and on
their respective frequencies $N_{1},\ldots,N_{K_n}$, recorded in the
sample $(X_1,\ldots,X_n)$. This will also serve as a tool for
predicting the value of the number of distinct species that will appear
$l$ times in the enlarged sample $(X_1,\ldots,X_{n+m})$, given the
observed sample $(X_1,\ldots,X_n)$.

Let $X_1^*,\ldots,X_{K_n}^*$ denote the labels identifying the $K_n$
distinct species detected in the sample $(X_1,\ldots,X_n)$. One can,
then, define
\[
L_m^{(n)}:=\sum_{i=1}^m \prod_{j=1}^{K_n}\mathbh{1}_{\{X_j^*\}
^c}(X_{n+i})=\operatorname{card}(\{X_{n+1},\ldots,X_{n+m}\}\cap
\{X_1^*,\ldots,X_{K_n}^*\}^c)
\]
as the number of observations from the additional sample of size $m$
that do not coincide with any of the $K_{n}$ distinct species in the
basic sample. Correspondingly, $X_{K_{n}+1}^*,\ldots
,X_{K_n+K_{m}^{(n)}}^*$ are the labels identifying the additional\vspace*{-2pt}
$K_{m}^{(n)}=K_{n+m}-K_n$ distinct species generated by these
$L_m^{(n)}$ observations. Then we can define
\[
S_{K_n+i}:=\sum_{j=1}^m\mathbh{1}_{\{X_{K_{n}+i}^*\}}(X_{n+j}),\qquad
S_q:=\sum_{j=1}^m\mathbh{1}_{\{X_q^*\}}(X_{n+j})
\]
for $i=1,\ldots,K_m^{(n)}$ and $q=1,\ldots,K_n$, where one obviously has
$\sum_{i=1}^{K_m^{(n)}}S_{K_n+i}=L_m^{(n)}$. For our purposes, it is
useful to resort to the decomposition $M_{l,n+m}=O_{l,m}+N_{l,m}$ where
\begin{equation}\label{eqoldnewsmallblocks}
O_{l,m}:=\sum_{q=1}^{K_n}\mathbh{1}_{\{l\}}(N_{q}+S_q),\qquad
  N_{l,m}:=\sum_{i=1}^{K_{m}^{(n)}}\mathbh{1}_{\{l\}}(S_{K_{n}+i})
\end{equation}
for any $l=1,\ldots,n+m$. It is apparent that $O_{l,m}=0$ for any
$l>n+m$ and $N_{l,m}=0$ for any $l>m$.
Hence, $O_{l,m}$ is the number of distinct species, among the $K_n$
detected in the basic sample $(X_1,\ldots,X_n)$, that have frequency
$l$ in the enlarged sample of size $n+m$. Analogously, $N_{l,m}$ is the
number of additional distinct species, generated by $L_m^{(n)}$
observations in $(X_{n+1},\ldots,X_{n+m})$, with frequency $l$ in the
enlarged sample. For notational convenience we introduce random
variables $O_{l,m}^{(n)}$ and $N_{l,m}^{(n)}$ that are defined in
distribution as follows:
\begin{eqnarray*}
\P\bigl[O_{l,m}^{(n)}=x\bigr]
&=&\P[O_{l,m}=x | K_n=j, \bm{N}=\bm{n}],\\
\P\bigl[N_{l,m}^{(n)}=y\bigr]
&=&\P[N_{l,m}=y | K_n=j, \bm{N}=\bm{n}]
\end{eqnarray*}
for any $1\le j\le n$, $\bm{n}\in\mathcal{D}_{n,j}$ and $n,m\ge1$.
Moreover, we set $\mathcal{C}_{j,r}$ as the space of all vectors
$\mathbf{c}^{(r)}=(c_{1},\ldots,c_{r})\in\{1,\ldots,j\}^{r}$ such that
$c_{i}\neq c_{\ell}$ for any $i\neq\ell$ and $\max_{1\le i\le
r}n_{c_i}\le l$. Finally,
\begin{eqnarray*}
&&I_{\sigma}\bigl(l,m,r,\bm{n},\mathbf{c}^{(r)}\bigr)\\
&&\qquad :=r!\pmatrix{m\vspace*{2pt}\cr l-n_{c_{1}},\ldots,l-n_{c_{r}},m-lr+ |\bm{n}_{\bm
{c}^{(r)}}|} \prod_{i=1}^{r}(n_{c_{i}}-\sigma)_{l-n_{c_{i}}},
\end{eqnarray*}
where $|\bm{n}_{\bm{c}^{(r)}}|:=\sum_{i=1}^r n_{c_i}$.
The next result provides an explicit expression for the $r$th factorial
moments of $O_{l,m}^{(n)}$ in terms of noncentral generalized factorial
coefficients defined by $\mathscr{C}(n,k;\gamma,\sigma):=(\sigma
t-\gamma)_{n}=\sum_{k=0}^{n}\mathscr{C}(n,k;\sigma,\gamma)(t)_{k}$ with
\mbox{$\sigma,\gamma\in\mathbb{R}$}. Recall also that $\mathscr
{C}(n,k;\sigma
,\gamma)=(1/k!)\sum_{j=0}^{k}(-1)^{j}{k\choose j}(-\sigma j-\gamma
)_{n}$ with the proviso $\mathscr{C}(0,0;\sigma,\gamma)=1$,
$\mathscr
{C}(n,0;\sigma)=(-\gamma)_{n}$ for any $n>0$ and
$\mathscr{C}(n,k;\sigma,\gamma)=0$ for any $k>n$.

\begin{thmm}\label{posteriormarginalold}
Let $(X_n)_{n\ge1}$ be an exchangeable sequence associated to a
Gibbs-type r.p.m. Then, for any $l=1,\ldots,n+m$, $r\geq1$ and $\bm
{n}\in\mathcal{D}_{n,j}$
%
\begin{eqnarray}\label{eqmarginalposteriorold2}
\E\bigl[\bigl(O_{l,m}^{(n)}\bigr)_{[r]}\bigr]&=&
\sum_{\mathbf{c}^{(r)}\in\mathcal{C}_{j,r}}
I_{\sigma}\bigl(l,m,r,\bm{n},\bm{c}^{(r)}\bigr)\nonumber\\
&&\hspace*{29pt}{}\times\sum_{k=0}^{m}\frac{V_{n+m,j+k}}{V_{n,j}} \\
&&\hspace*{56pt}{}\times\frac{\mathscr
{C}
(m-rl+|\bm{n}_{\bm{c}^{(r)}}|,k;\sigma,-n+|\bm{n}_{\bm{c}^{(r)}}|
+(j-r)\sigma)}{\sigma^{k}}.\nonumber
\end{eqnarray}
%
\end{thmm}

It is worth observing that the moments in \eqref
{eqmarginalposteriorold2}, for any $r\ge1$, characterize the
distribution of $O_{l,m}^{(n)}$. Such a distribution is interpretable
as the posterior probability distribution, given the observations
$(X_1,\ldots,X_n)$, of the number of distinct species that (i) appear
with frequency $l$ in a sample of size $n+m$; (ii) had been already
detected within $(X_1,\ldots,X_n)$. Therefore, we will refer to
$O_{l,m}^{(n)}$ as the number of ``old'' species with frequency $l$.\vadjust{\goodbreak}
The Bayesian nonparametric estimator, under a quadratic loss function,
coincides with the expected value of $O_{l,m}^{(n)}$ and is easily
recovered from \eqref{eqmarginalposteriorold2}.
%
\begin{cor}\label{estimatorold}
Let $(X_n)_{n\ge1}$ be an exchangeable sequence associated to a
Gibbs-type r.p.m. Conditionally on a sample $(X_1,\ldots,X_n)$, the
expected number of ``old'' distinct species that appear with frequency
$l$, for any $l=1,\ldots,n+m$, in a sample of size $n+m$ is given by\vspace*{-1pt}
\begin{eqnarray}\label{eqbayesestimatorold}
\hat{O}_{l,m}^{(n)}&:=&\E\bigl[O_{l,m}^{(n)}\bigr]\nonumber\\
&=&\sum_{t=1}^{l}\pmatrix{m\vspace*{2pt}\cr
l-t}m_{t}(t-\sigma)_{l-t}
\\
&&\hspace*{16pt}{}\times
\sum_{k=0}^{m}\frac{V_{n+m,j+k}}{V_{n,j}}
\frac{\mathscr{C}(m-(l-t),k;\sigma,-n+t+(j-1)\sigma
)}{\sigma^{k}},\nonumber
\end{eqnarray}
with $m_{t}\geq0$ being the number of distinct species with frequency
$t$ observed in the basic sample, namely, $m_t=\sum_{i=1}^{K_n}\indic
_{\{t\}}(N_{i})$. Moreover, $(K_{n},M_{1,n},\ldots,M_{l,n})$ is
sufficient for predicting $O_{l,m}^{(n)}$ over the whole sample of size $n+m$.
\end{cor}

An analogous result of Theorem~\ref{posteriormarginalold} can be
established for $N_{l,m}^{(n)}$. Indeed, if we set
\[
J_{\sigma}(l,m,r):=\pmatrix{m\vspace*{2pt}\cr l,\ldots,l,m-rl} [(1-\sigma
)_{l-1}]^r,
\]
one can show the following theorem.

\begin{thmm}\label{posteriormarginalnew}
Let $(X_n)_{n\ge1}$ be an exchangeable sequence associated to a
Gibbs-type r.p.m. Then, for any $l=1,\ldots,m$ and $r\geq1$,
%
%
\begin{equation}
\label{eqmarginalposteriornew2}
\quad\E\bigl[\bigl(N_{l,m}^{(n)}\bigr)_{[r]}\bigr]=
J_{\sigma}(l,m,r)\sum_{k=0}^{m-rl}\frac{V_{n+m,j+k+r}}{V_{n,j}}
\frac
{\mathscr{C}(m-rl,k;\sigma,-n+j\sigma)}{\sigma^{k}}.
\end{equation}
%
\end{thmm}

Hence, \eqref{eqmarginalposteriornew2} characterizes the
probability distribution of $N_{l,m}^{(n)}$. This can be seen as the
posterior probability distribution, conditional on the observations
$(X_1,\ldots,X_n)$, of the number of distinct species that (i) appear
with frequency~$l$ in a sample of size $n+m$; (ii) do not coincide with
any of the $K_n$ distinct species already detected within $(X_1,\ldots
,X_n)$. For this reason $N_{l,m}^{(n)}$ is referred to as the number of
``new'' species with frequency $l$. Thus, the Bayesian nonparametric
estimator, under a quadratic loss function, is easily recovered from~\eqref{eqmarginalposteriornew2}.

\begin{cor}\label{estimatornew}
Let $(X_n)_{n\ge1}$ be an exchangeable sequence associated to a
Gibbs-type r.p.m. Conditionally on a sample $(X_1,\ldots,X_n)$, the\vadjust{\goodbreak}
expected number of ``new'' distinct species that appear with frequency
$l$, for any $l=1,\ldots,m$, in a sample of size $n+m$ is given by
%
%
\begin{eqnarray}
\label{eqbayesestimatornew}
\hat{N}_{l,m}^{(n)}&:=&\E\bigl[N_{l,m}^{(n)}\bigr]
\nonumber
\\[-8pt]
\\[-8pt]
\nonumber
&\phantom{:}=&
\pmatrix{m\vspace*{2pt}\cr l}(1-\sigma)_{l-1}
\sum_{k=0}^{l}\frac{V_{n+m,j+k+1}}{V_{n,j}}
 \frac{\mathscr{C}(m-l,k;\sigma,-n+j\sigma)}{\sigma^{k}}.
\end{eqnarray}
Hence, $K_{n}$ is sufficient for predicting $N_{l,m}^{(n)}$.
\end{cor}

\begin{rmk}
According to the definition of the random variable $N_{l,m}^{(n)}$, one has
%
%
\begin{equation}\label{eqbayesestimator}
\hat{E}^{(n)}_{m}:=\E\bigl[K_{m}^{(n)} |
K_{n}=j,\bm{N}=\bm{n}\bigr]=\sum_{l=1}^{m}\hat{N}_{l,m}^{(n)},
\end{equation}
providing an alternative derivation of the Bayesian nonparametric
estimator for the number of ``new'' distinct species obtained in \cite
{Lij308}. A detailed discussion of the estimator \eqref
{eqbayesestimator} and its relevance in genomics can be found in~\cite
{lmp07bmc}.\looseness=-1
\end{rmk}

At this point we turn our attention to characterizing the following
random variable:
%
%
\begin{equation}\label{eqdefoldnew}
M^{(n)}_{l,m}\stackrel{\mathrm{d}}{=}O_{l,m}^{(n)}+N_{l,m}^{(n)},
\end{equation}
whose probability distribution coincides with the distribution of the
number $M_{l,n+m}$ of clusters of size $l$ featured by $(X_1,\ldots
,X_{n+m})$ conditional on the basic sample
$(X_1,\ldots,X_n)$. In particular, if we set
\begin{eqnarray*}
&&H_{\sigma}\bigl(l,m,r,t,\bm{n},\bm{c}^{(t)}\bigr)\\
&&\qquad :=t!\pmatrix{m\vspace*{2pt}\cr l,\ldots,l,l-n_{c_{1}},\ldots
,l-n_{c_{t}},m-rl+\bigl|\bm
{n}_{\bm{c}^{(t)} }\bigr|}\\
&&\quad\qquad{}
\times[(1-\sigma)_{l-1}]^{r-t}\prod_{i=1}^{t}(n_{c_{i}}-\sigma
)_{l-n_{c_{i}}},
\end{eqnarray*}
an analogous result of Theorems~\ref{posteriormarginalold} and~\ref{posteriormarginalnew} can be established for $M^{(n)}_{l,m}$.

\begin{thmm}\label{posteriormarginaloldnew}
Let $(X_n)_{n\ge1}$ be an exchangeable sequence associated to a
Gibbs-type r.p.m. Then, for any $l=1,\ldots,m+n$ and $r\geq1$,
\begin{eqnarray}
\label{eqmarginalposterioroldnew2}
&&\E\bigl[\bigl(M_{l,m}^{(n)}\bigr)_{[r]}\bigr]\nonumber\\
&&\qquad=\sum
_{t=0}^{r}\pmatrix{r\vspace*{2pt}\cr t}\sum_{\mathbf{c}^{(t)}\in\mathcal
{C}_{j,t}}H_{\sigma}\bigl(l,m,r,t,\bm{n},\bm{c}^{(t)}\bigr)
\nonumber
\\[-8pt]
\\[-8pt]
\nonumber
&&\hspace*{68pt}\qquad\quad{} \times\sum_{k=0}^{m-rl+|\bm{n}_{\bm{c}^{(t)} }|}\frac
{V_{n+m,j+k+r-t}}{V_{n,j}}\\
&&\hspace*{79pt}\qquad{}\times\frac{\mathscr{C}(m-rl+|\bm{n}_{\bm{c}^{(t)}
}|,k;\sigma, -n+|\bm{n}_{\bm{c}^{(t)} }|+(j-t)\sigma)
}{\sigma^{k}}.\nonumber
\end{eqnarray}
\end{thmm}

Hence, \eqref{eqmarginalposterioroldnew2} characterizes the
probability distribution of $M^{(n)}_{l,m}$. This is interpreted as the posterior probability
distribution, given the observation $(X_{1},\ldots,X_{n})$, of the
number of distinct species that appear with frequency $l$ in a sample
of size $n+m$. Thus, the Bayesian nonparametric estimator,
under a quadratic loss function, is easily recovered from \eqref
{eqmarginalposterioroldnew2}. Clearly, according to \eqref{eqdefoldnew},
this also corresponds to the sum of the estimators in \eqref
{eqbayesestimatorold} and \eqref{eqbayesestimatornew}.


\section{Illustrations}\label{sec3} We now apply the general results of Section~\ref{sec2}
and specialize them to some noteworthy examples of Gibbs-type models.
We will devote particular attention to the two-parameter
Poisson--Dirichlet process since it is particularly suited for species
sampling applications in general \cite{Lij207} and for genomic
applications in particular, as will be seen in Section~\ref{sec4}.

\subsection{The Dirichlet process}\label{sec3.1} Denote the signless Stirling number
of the first kind by $|s(n,k)|$ and recall that
$
\lim_{\sigma\to0}\sigma^{-k} \Ccr(n,k,\sigma)=|s(n,k)|
$
for any $n\ge1$ and $1\le k\le n$. Now, let $\tilde P$ be a Dirichlet
process with parameter $\theta$ and, considering the form of the
$V_{n,j}$ weights and Theorem~\ref{priormarginal}, one readily obtains
\[
\E\bigl[(M_{l,n})_{[r]}\bigr]=\frac{(n)_{[rl]}}{l^r (\theta)_n}
\sum_{j=1}^{n-rl+r}
\theta^j |s(n-rl,j-r)|=\frac{(n)_{[rl]}}{l^r (\theta)_n} (\theta
)_{[n-rl]}.
\]
Using the classical sieve formula, one easily shows the following,
which appears to be new even in the case of Ewens partitions with the
exception of the case $l=1$ obtained in \cite{Ewe98}.

\begin{prp}\label{prpdirprior}
If $(X_n)_{n\ge1}$ is an exchangeable sequence associated to a
Dirichlet process with parameter $\theta>0$, then, for any $n \geq1$
and $l=1, \ldots, n$, the distribution of $M_{l,n}$ is of the form
%
%
\begin{equation}
\P[M_{l,n}=m_{l}]
=\frac{n!}{m_l!(\theta)_n}\frac{\theta^{m_l}}{l^{m_l}} \sum
_{t=0}^{[n/l]-m_l}
\frac{(-1)^t(\theta)_{[n-m_l-tl]}}{(n-m_ll-tl)!}\frac{\theta^t}{l^t}.
\label{eqmarginaldir}
\end{equation}
\end{prp}

On the basis of the result stated in Proposition~\ref{prpdirprior},
one can derive the asymptotic behavior of $M_{l,n}$, namely, that, for
any $l\geq1$,
%
%
\begin{equation}\label{limit2}
M_{l,n}\stackrel{\mathrm{d}}{\longrightarrow} W_l
\end{equation}
as $n\rightarrow+\infty$, where $W_l$ is a random variable distributed
according to a Poisson distribution with parameter $\theta/l$. The
limit result
\eqref{limit2} is known in the literature and has been originally
obtained in \cite{Arr92,Bar92}. See also \cite{Arr03} and
references therein.

Turning attention to the conditional case, one can easily derive the
following results. Theorem~\ref{posteriormarginalold} provides an
expression for the probability distribution of $O^{(n)}_{l,m}$, that is,
\begin{eqnarray*}
&&\P\bigl[O^{(n)}_{l,m}=m_{l}\bigr]\\
&&\qquad =\sum_{t=0}^{m-m_{l}} (-1)^{t} \pmatrix{m_{l}+t\vspace*{2pt}\cr t } \sum_{\bm
{c}^{(m_{l}+t)}\in\mathcal{C}_{j,m_{l}+t}}\frac{m!}{\prod
_{i=1}^{m_l+t}(l-n_{c_i})!
(m-\nu_t)!}\\[7pt]
&&\hspace*{20pt}\qquad\quad{}\times\prod_{i=1}^{m_{l}+t}(n_{c_{i}})_{l-n_{c_{i}}}\frac{
(\theta
+n-\sum_{i=1}^{m_{l}+t}n_{c_{i}})_{m-\nu_t}}{(\theta+n)_{m}},
\end{eqnarray*}
where we set $\nu_t=\sum_{i=1}^{m_{l}+t}(l-n_{c_{i}})$. Analogously,
Theorem~\ref{posteriormarginalnew} provides an
expression for the probability distribution of $N^{(n)}_{l,m}$, that is,
\[
\P\bigl[N^{(n)}_{l,m}=m_{l}\bigr]=\frac{\theta
^{m_{l}}}{t^{m_l}}\sum
_{t=0}^{m-m_{l}}
\biggl(-\frac{\theta}{l}\biggr)^{t} \frac{m!}{t!m_l!(m-lm_l-lt)!}
\frac{(\theta+n)_{m-lm_{l}-lt}}{(\theta+n)_{m}}.
\]
Similarly, according to Corollaries~\ref{estimatorold} and~\ref
{estimatornew}, and using the limiting result for
noncentral generalized factorial coefficients
\[
\lim_{\sigma\rightarrow0}\frac{\mathscr{C}(n,k;\sigma,\gamma
)}{\sigma
^{k}}=\sum_{i=k}^{n}\pmatrix{n\vspace*{2pt}\cr i}|s(i,k)|(-\gamma)_{n-i},
\]
the Bayesian estimators of the number of ``old'' and of ``new'' species
of size $l$ generated by $(X_1,\ldots,X_{n+m})$, conditional on
$(X_1,\ldots,X_n)$, are
given by
\begin{equation}\label{eqestimold}
\hat{O}_{l,m}^{(n)}=\sum_{t=1}^{l}\pmatrix{m\vspace*{2pt}\cr l-t}m_{t}(t)_{l-t}\frac
{(\theta+n-t)_{m-(l-t) }}{(\theta+n)_{m }}
\end{equation}
and
\begin{equation}\label{eqestimnew}
\hat{N}_{l,m}^{(n)}=(l-1)!\pmatrix{m\vspace*{2pt}\cr l}\frac{\theta}{(\theta+n+m-l)_{l}}.
\end{equation}
In particular, from \eqref{eqestimold} and \eqref{eqestimnew} the
Bayesian estimator of the number of clusters of size $l$ over an
enlarged sample of size $n+m$, conditional on the partition structure
of the $n$ observed data, is given in the following proposition.

\begin{prp}\label{prpdir2}
If $(X_n)_{n\ge1}$ is an exchangeable sequence associated to a
Dirichlet process with parameter $\theta$, then
\[
\hat{M}_{l,m}^{(n)}=\pmatrix{{m}\vspace*{2pt}\cr{l}} \frac{\theta (l-1)!}{(\theta+n+m-l)_l}
+\sum_{t=1}^l\pmatrix{{m}\vspace*{2pt}\cr{l-t} } m_t (t)_{l-t} \frac{(\theta
+n-t)_{m-l+t}}{(\theta+n)_m}
\]
for any $l\in\{1,\ldots,n+m\}$.
\end{prp}

Finally, by combining \eqref{eqmarginalposteriornew2} and \eqref
{eqmarginalposterioroldnew2} a simple limiting argument leads to
show that, as $m\rightarrow+\infty$ and for any $l\ge1$,
$N_{l,m}^{(n)}\stackrel{\mathrm{d}}{\longrightarrow} W^{(n)}_l$
and
%
%
\begin{equation}\label{eqlimdirpost}
M_{l,m}^{(n)}\stackrel{\mathrm{d}}{\longrightarrow} W^{(n)}_l,
\end{equation}
where $W^{(n)}_l$ is a random variable distributed according to a
Poisson distribution with parameter $(\theta+n)/l$. Clearly, \eqref
{eqlimdirpost} reduces to \eqref{limit2} in the unconditional case
corresponding to $n=0$.

\subsection{The two-parameter Poisson--Dirichlet process}\label{sec3.2} The Pitman
model with parameters $(\sigma,\theta)$ in \eqref{eqeppf}, or
$\operatorname{PD}(\sigma,\theta)$ process, stands out for its analytical tractability
and for its modeling flexibility. In particular, within the species
sampling context, the presence of the additional parameter $\sigma\in
(0,1)$, w.r.t.~the simple Dirichlet model, allows to model more
effectively both the clustering structure featured by the $X_i$'s and
the growth rate of $K_n$. Therefore, given its importance, we devote
special attention to this process. A few additional asymptotic results
that complement, for the specific case we are analyzing, those recalled
in Section~\ref{sec2} for general Gibbs-type r.p.m.'s are of particular interest.

\subsubsection{Distributional results}\label{sec3.2.1} Let us first state a result
concerning the unconditional distribution of $M_{l,n}$, namely, the
number of clusters with frequency $l$ in a sample of size $n$.

\begin{prp} \label{prp2pd0} Let $(X_n)_{n\ge1}$ be an exchangeable
sequence associated to a $\operatorname{PD}(\sigma,\theta)$ process with $\sigma\in
(0,1)$ and $\theta>-\sigma$. Then,
\begin{eqnarray}
\label{eqmarginalpd}
\P[M_{l,n}=m_{l}]
&=&\sum_{t=0}^{n-m_{l}}(-1)^{t}\frac{n!}{t!m_{l}!(n-lm_{l}-lt)!} \sigma
^{m_l+t}\biggl(\frac{\theta}{\sigma}\biggr)_{m_{l}+t}
\nonumber
\\[-8pt]
\\[-8pt]
\nonumber
&&\hspace*{18pt}{}\times \biggl(\frac{(1-\sigma)_{l-1}}{l!}\biggr)^{m_{l}+t}
\frac{(\theta+(m_{l}+t)\sigma)_{n-lm_{l}-lt} }{(\theta)_{n}}.
\end{eqnarray}
\end{prp}

Hence, \eqref{eqmarginalpd} provides the marginal distribution of the
Pitman sampling formula~\eqref{eqgibbssampling}, corresponding to
$V_{n,j}=\sigma^j(\theta/\sigma)_j/(\theta)_n$, and, to the authors'
knowledge, it is not explicitly reported in the literature.

Turning attention to the conditional case, one can easily derive the
following results.

\begin{prp}\label{prppdold}
Let $(X_n)_{n\ge1}$ be an exchangeable sequence associated to a
$\operatorname{PD}(\sigma,\theta)$ process with $\sigma\in(0,1)$ and $\theta
>-\sigma$. Then,
\begin{eqnarray}
\label{eqold2pd}
&&\P\bigl[O^{(n)}_{l,m}=m_{l}\bigr]\nonumber\\
&&\qquad=\sum_{t=0}^{m-m_{l}}(-1)^{t}
\pmatrix{m_{l}+t\vspace*{2pt}\cr t}
\nonumber
\\
&&\hspace*{19pt}\qquad\quad{}\times \sum_{\bm{c}^{(m_{l}+t)}\in\mathcal
{C}_{j,m_{l}+t}}\pmatrix{m\vspace*{2pt}\cr
l-n_{c_{1}},\ldots,l-n_{c_{m_{l}+t}},\sum_{i=1}^{m_{l}+t}(l-n_{c_{i}})
}\\
&&\hspace*{19pt}\qquad\quad{}\times
\prod_{i=1}^{m_{l}+t}(n_{c_{i}}-\sigma)_{l-n_{c_{i}}} \nonumber\\
&&\hspace*{19pt}\quad\qquad{}\times\frac{(\theta+n-\sum
_{i=1}^{m_{l}+t}n_{c_{i}}+(m_{l}+t)\sigma)_{m-\sum
_{i=1}^{m_{l}+t}(l-n_{c_{i}})}}{(\theta+n)_{m}}\nonumber
\end{eqnarray}
for any $l\in\{1,\ldots,n\}$ and $m_l\in\{1,\ldots,n\}$ such that
$m_ll\le n$.
\end{prp}

From \eqref{eqold2pd} one can deduce a completely explicit expression
for the Bayesian estimator of the number of ``old'' species with
frequency $l$ in the whole sample $X_1,\ldots,X_{n+m}$, namely,
%
%
\begin{equation}\label{bayesestimatorfav1}
\hat{O}_{l,m}^{(n)}=\E\bigl[O_{l,m}^{(n)}\bigr]=
\sum_{t=1}^{l}\pmatrix{m\vspace*{2pt}\cr l-t}m_{t}(t-\sigma)_{l-t}\frac{(\theta
+n-t+\sigma)_{m-(l-t) }}{(\theta+n)_{m }},
\end{equation}
which can be readily used in applications, as will be shown in Section~\ref{sec4}.
In a~similar fashion it is possible to deduce the distribution of the
number of ``new'' species that will appear $l$ times in
$(X_{n+1},\ldots
,X_{n+m})$ conditional on the observations $(X_1,\ldots,X_n)$. Indeed,
one can show the following:

\begin{prp}\label{prp2pd2}
Let $(X_n)_{n\ge1}$ be an exchangeable sequence associated to a
$\operatorname{PD}(\sigma,\theta)$ process with $\sigma\in(0,1)$ and $\theta
>-\sigma$. Then,
\begin{eqnarray}\label{eqdistributionnewuniqueness}
\P\bigl[N_{l,m}^{(n)}=m_{l}\bigr]
&=&\sum_{t=0}^{m-m_{l}}(-1)^{t}\pmatrix{m\vspace*{2pt}\cr t,m_{l},m-lm_{l}-lt }\prod
_{i=0}^{m_{l}+t-1}(\theta+j\sigma+i\sigma)
\nonumber
\\[-8pt]
\\[-8pt]
\nonumber
&&\hspace*{20pt}{}\times\biggl(\frac{(1-\sigma)_{l-1}}{l!}\biggr)^{m_{l}+t}
\frac{(\theta+n+(m_{l}+t)\sigma)_{m-l(m_{l}+t)}}{(\theta+n)_{m}},
\end{eqnarray}
for any $n\geq1$, $j=1,\ldots,n$, $l\ge1$ and $m_l\ge1$ such that
$m_l l\le m$.
\end{prp}

\begin{rmk}
One can alternatively prove \eqref{eqdistributionnewuniqueness} by
relying on the so-called quasi-conjugacy property of the
two-parameter Poisson--Dirichlet process, a concept introduced in
\cite
{Lij308}. Indeed, it suffices to marginalize an \textit{updated}
Pitman sampling formula and \eqref{eqdistributionnewuniqueness}
easily follows. Moreover, if $n=j=0$ in \eqref
{eqdistributionnewuniqueness}, one recovers the marginal
distribution of $M_{l,n}$ as described in \eqref{eqmarginalpd} and,
if one additionally sets $\sigma=0$, the distribution of $M_{l,n}$
corresponding to the Ewens partition in \eqref{eqmarginaldir} is obtained.
\end{rmk}

The Bayesian estimator for the number of ``new'' species with frequency
$l$ over the enlarged sample $n+m$ coincides with
%
%
\begin{equation}
\label{bayesestimatorfav2}
\hat{N}_{l,m}^{(n)}=\E\bigl[N_{l,m}^{(n)}\bigr]=\pmatrix{m\vspace*{2pt}\cr l}(1-\sigma)_{l-1
}(\theta+j\sigma)\frac{(\theta+n+\sigma)_{m-l }}{(\theta+n)_{m}}
\end{equation}
for any $l\in\{1,\ldots,m\}$. Having determined $\hat{O}_{l,m}^{(n)}$
and $\hat{N}_{l,m}^{(n)}$, one finds out that a Bayesian estimator of
the total number of species with frequency $l$ among $(X_1,\ldots
,X_{n+m})$, given $(X_1,\ldots,X_n)$, is given by the following:

\begin{prp}\label{prp2pd3}
If $(X_n)_{n\ge1}$ is an exchangeable sequence with $\tilde P$ in
\eqref{eqexchangeabl} being the $\operatorname{PD}(\sigma,\theta)$ process, for any
$l=1,\ldots,n+m$,
\begin{eqnarray}\label{pippo}
\hat{M}_{l,m}^{(n)}
&=&
\sum_{t=1}^{l}\pmatrix{m\vspace*{2pt}\cr l-t}m_{t}(i-\sigma)_{l-t}\frac{(\theta
+n-t+\sigma)_{m-(l-t) }}{(\theta+n)_{m }}
\nonumber
\\[-8pt]
\\[-8pt]
\nonumber
&&{}  +\pmatrix{m\vspace*{2pt}\cr l}(1-\sigma)_{l-1 }(\theta+j\sigma
)\frac
{(\theta+n+\sigma)_{m-l }}{(\theta+n)_{m}}.
\end{eqnarray}
\end{prp}

Of course, Theorem~\ref{posteriormarginaloldnew} allows a direct
evaluation of $\hat{M}_{l,m}^{(n)}$ above and yields moments of any
order $r\geq1$
of $M_{l,m}^{(n)}$.

\subsubsection{Asymptotics}\label{sec3.2.2} We now study the asymptotic behavior of
$M_{l,m}^{(n)}$ and $N_{l,m}^{(n)}$, as $m\to\infty$. However, before
proceeding, let us first recall a well-known result concerning the
asymptotics of $M_{l,n}$ as $n$ increases. To this end, let $f_{\sigma
}$ be the density function of a positive $\sigma$-stable random
variable and $Y_{q}$, for any $q\geq0$, a positive random variable with
density function
\[
f_{Y_{q}}(y)=\frac{\Gamma(q\sigma+1)}{\sigma\Gamma
(q+1)}y^{q-1/\sigma
-1}f_{\sigma}(y^{-1/\sigma}).
\]
Then, for any $l\geq1$
\[
\frac{M_{l,n}}{n^{\sigma}}\stackrel{\mathrm{d}}{\longrightarrow
}\frac
{\sigma(1-\sigma)_{(l-1)}}{l!}Y_{\theta/\sigma}
\]
as $n\rightarrow+\infty$. See \cite{Pit06} for details. We now
provide a new result concerning the limiting behavior in the
conditional case and, specifically, of\vadjust{\goodbreak} $M_{l,m}^{(n)}$ and of
$N_{l,m}^{(n)}$ as $m\to\infty$. It will be shown that they converge in
distribution to the same random element that still depends on $Y_q$ for
a suitable choice of $q$.

\begin{thmm}\label{limitbayesestimatoruniqueness}
Let $(X_n)_{n\ge1}$ be an exchangeable sequence associated to a
$\operatorname{PD}(\sigma,\theta)$ process. For any $1\le j\le n$ and $l\ge1$, one has
%
%
\begin{equation}
\label{limituniqueness1new}
\frac{N_{l,m}^{(n)}}{m^{\sigma}}
\stackrel{d}{\longrightarrow}
\frac{\sigma(1-\sigma)_{l-1}}{l!}Z_{n,j}
\end{equation}
as $m\rightarrow+\infty$, where $Z_{n,j}\stackrel{d}{=}B_{j+\theta
/\sigma,n/\sigma-j}Y_{(\theta+n)/\sigma}$ and $B_{j+\theta/\sigma
,n/\sigma-j}$ is a beta random variable with parameters $(j+\theta
/\sigma,n/\sigma-j)$ independent of $Y_{(\theta+n)/\sigma}$. Moreover,
%
%
\begin{equation}\label{limituniqueness1}
\frac{M_{l,m}^{(n)}}{m^{\sigma}}\stackrel{d}{\longrightarrow}
\frac{\sigma(1-\sigma)_{l-1}}{l!}Z_{n,j}
\end{equation}
as $m\rightarrow+\infty$.
\end{thmm}

The limit in \eqref{limituniqueness1new} and \eqref
{limituniqueness1} implies that $K_{n}$ is asymptotically sufficient
for predicting the conditional number of
distinct species with frequency $l$ to be generated by the additional
sample $(X_{n+1},\ldots,X_{n+m})$ as its size $m$ increases. Such a
limit involves the beta-tilted
random variable $Z_{n,j}$, originally introduced in \cite{Fav09} by
investigating the asymptotic behavior of the conditional number of
``new'' distinct species $K_{m}^{(n)}$
generated by the additional sample as its size $m$ increases.
Specifically,
\[
\frac{K_{m}^{(n)}}{m^{\sigma}}\stackrel{}{\rightarrow}Z_{n,j},
\]
almost surely, as $m\rightarrow+\infty$. It is worth noting that
beta-tilted random variables of similar type have been recently the
object of a thorough investigation in \cite{Jam10} in the context of
the so-called Lamperti-type laws.

\begin{rmk}
Note that from \eqref{limituniqueness1new} and \eqref
{limituniqueness1} one obtains the unconditional result of \cite
{Pit06} by setting $n=j=0$. Moreover, one recovers a result in \cite
{Fav09}, which states that, conditional on $(X_1,\ldots,X_n)$,
$m^{-\sigma} K_m^{(n)}\stackrel{\mathrm{d}}{\to} Z_{n,j}$, as $m\to\infty$.
Indeed, $K_m^{(n)}=\sum_{l=1}^{L_m^{(n)}}N_{l,m}^{(n)}$ and $L_m^{(n)}$
diverges as $m \to+ \infty$: hence, the limit in distribution for
$K_m^{(n)}$ can be deduced from \eqref{limituniqueness1new} upon
noting that $\sum_{l\ge1} (l!)^{-1} \sigma(1-\sigma)_{l-1}=1$.
\end{rmk}

\subsection{The Gnedin model}\label{sec3.3} Consider now the Gnedin model \eqref
{eqgnedinweight} with parameters $\zeta=0$ and $\gamma\in[0,1)$. The
corresponding random partition is representable as a mixture partitions
of the type \eqref{eqeppf}, however, with parameters $(-1,\kappa)$,
each of which generates a partition with a finite number of blocks
$\kappa$. The mixing distribution for the total number of blocks is
$p(\kappa)=\indic_{\{1,2,\ldots\}}(\kappa)\gamma (1-\gamma)_{\kappa-1}/\kappa!$.\vadjust{\goodbreak}

\begin{prp}\label{prppriorgnedin}
Let $(X_n)_{n\ge1}$ be an exchangeable sequence associated to the
Gnedin model with parameters $(0,\gamma)$. Then
\begin{eqnarray}\label{eqpriormomgnedin}
&&\E\bigl[(M_{l,n})_{[r]}\bigr]\nonumber\\
&&\qquad=\indic_{\{rl\}}(n)
\frac
{r! l  (\gamma)_{rl-r}(1-\gamma)_{r-1}}{(1+\gamma)_{rl-1}}
+\indic_{\{rl+1,\ldots\}}(n) \frac{n(\gamma)_{rl-r}(1-\gamma
)_{r}}{(1+\gamma)_{n-1}}\\
&&\qquad\quad{}\times
\sum_{k=0}^{n-rl-1}\pmatrix{{n-rl-1}\vspace*{2pt}\cr{k}}
\frac{(r+k)!}{(1+k)!}
(\gamma+rl-r)_{n-rl-1-k}(r+1-\gamma)_k.\nonumber
\end{eqnarray}
\end{prp}

From \eqref{eqpriormomgnedin} one can determine the probability
distribution of $M_{l,n}$. Indeed, if $n/l\notin\N$, then
\begin{eqnarray*}
&&\hspace*{-4pt}\P[M_{l,n}=m_{l}]\\
&&\hspace*{18pt}=\frac{\indic_{\{1,\ldots,n\}}(lm_l)
n}{m_l!(1+\gamma)_{n-1}}\\
&&\hspace*{27pt}{}\times \sum_{r=m_l}^{[n/l]}\frac{(-1)^{r-m_l}}{(r-m_l)!} (\gamma
)_{rl+r}(1-\gamma)_r\\
&&\hspace*{58pt}{}\times \sum_{k=0}^{n-rl-1}\pmatrix{{n-rl-1}\vspace*{2pt}\cr{k} }\frac
{(r+k)!}{(1+k)!}(\gamma+rl+r)_{n-rl-1-k}
(r+1-\gamma)_{k}.
\end{eqnarray*}
On the other hand, if $n/l\in\N$, then
\begin{eqnarray*}
\P[M_{l,n}=m_{l}]&=&\frac{\indic_{\{1,\ldots,n\}}(lm_l)
n}{m_l!(1+\gamma
)_{n-1}}\\
&&{}\times \Biggl \{
(-1)^{{n}/{l}-m_l} \frac{({n}/{l}-1)!(\gamma
)_{n-{n}/{l}}
(1-\gamma)_{{n}/{l}-1}}{({n}/{l}-m_l)!}
\\
&&\hspace*{19pt}{}+
 \sum_{r=m_l}^{{n}/{l}-1}\frac{(-1)^{r-m_l}}{(r-m_l)!} (\gamma
)_{rl+r}(1-\gamma)_r\\
&&\hspace*{60pt}{}\times \sum_{k=0}^{n-rl-1}\pmatrix{{n-rl-1}\vspace*{2pt}\cr{k}} \frac{(r+k)!}{(1+k)!}
\\
&&\hspace*{100pt}{}\times(\gamma+rl+r)_{n-rl-1-k}(r+1-\gamma)_{k}\Biggr\}.
\end{eqnarray*}


Moreover, for any $l\geq1$,
%
%
\begin{equation}\label{limit3}
M_{l,n}\stackrel{\mathrm{d}}{\longrightarrow} 0
\end{equation}
as $n\rightarrow+\infty$. Note that the limiting result in \eqref
{limit3} is not surprising since a Gnedin r.p.m.~induces a random
partition of $\N$ into an almost surely finite number of blocks even
though with infinite expectation \cite{Gne10}.

As for the posterior distribution of the number of clusters of size
$l$, we now use the general results outlined in Section~\ref{sec2} to provide
some explicit forms for the distribution of $O_{l,m}^{(n)}$ and $N_{l,m}^{(n)}$.

\begin{prp}\label{prpgnedin} Let $(X_n)_{n\ge1}$ be an exchangeable
sequence associated to the Gnedin model with parameters $\zeta=0$ and
$\gamma\in[0,1)$. Then,
\begin{eqnarray*}
\P\bigl[O^{(n)}_{l,m}=m_{l}\bigr]&=&\frac{\indic_{\{1,\ldots,n\}}(lm_l)  m!}{(n)_m(\gamma
+n)_m} \\
&&{}\times \sum_{r=m_l}^{[n/l]} (-1)^{r-m_l}\pmatrix{{r}\vspace*{2pt}\cr{m_l}}(m+n+j-r-rl-1)!
\\
&&\hspace*{8pt}\qquad{}\times \sum_{\bm{c}^{(r)}\in\mathcal{C}_{j,r}} \frac{1}{(m-rl+|\bm
{n}_{\bm
{c}^{(r)}}|)!} \prod_{i=1}^r
\frac{(n_{c_i}-\sigma)_{l-n_{c_i}}}{(l-n_{c_i})!}\\
&&\hspace*{73pt}{}\times \sum_{k=0}^{m-rl+|\bm{n}_{\bm{c}^{(r)}}|}\pmatrix{{m-rl+|\bm
{n}_{\bm{c}^{(r)}}|}\vspace*{2pt}\cr{k}}
\\
&&\hspace*{113pt}\qquad{}\times  \frac{(j)_k  (\gamma+n-j)_{m-k} }{(n-|\bm{n}_{\bm{c}^{(r)}}|+j-r-1+k)!}.
\end{eqnarray*}
Moreover,
\begin{eqnarray*}
\P\bigl[N^{(n)}_{l,m}=m_{l}\bigr]
&=&\frac{\indic_{\{1,\ldots,m\}}(lm_l)  m!}{(n)_m(\gamma+n)_m} \\
&&{}\times \sum
_{r=m_l}^{[n/l]}
(-1)^{r-m_l} \frac{(m-rl+n+j)!}{(r-m_l)!(m-rl)!}\\
&&\hspace*{30pt}{}\times  \sum_{k=0}^{m-rl} \pmatrix{{m-rl}\vspace*{2pt}\cr{k}}  \frac{(\gamma
+n-j)_{m-r-k}(j)_{k+r}(j-\gamma)_{k+r}}{(n+j+k)!}.
\end{eqnarray*}
\end{prp}

One can further deduce the conditional expected values of
$O_{l,m}^{(n)}$ and of $N_{l,m}^{(n)}$ which take on the following forms:
\begin{eqnarray*}
\hat{O}_{l,m}^{(n)}
&=&\frac{1}{(n)_{m}(\gamma+n)_m}\\
&&{}\times \sum_{t=1}^{l}m_t  \pmatrix{m \vspace*{2pt}\cr l-t} (t+1)_{l-t}  (m+n+j-l-2)!\\
&&\hspace*{25pt}{} \times \sum_{k=0}^{m-l+t}\pmatrix{{m-l+t}\vspace*{2pt}\cr{k}} \frac{(\gamma
+n-j)_{m-k}(j)_k}{(n+j-t-2+k)!},
\\
\hat{N}_{l,m}^{(n)}
&=&\frac{m!(1+\gamma)_{n-1} (n+j)_{m-l}}{(n)_m (\gamma+n)_m}
\sum_{k=0}^{m-l}\pmatrix{{m-l}\vspace*{2pt}\cr{k}} \frac{(j)_k(j-\gamma)_{k+1}}{(n+j)_k}.
\end{eqnarray*}
As in previous examples, these quantities can, then, be used in order
to provide a Bayesian estimator $\hat{M}_{l,m}^{(n)}=\hat
{O}_{l,m}^{(n)}+\hat{N}_{l,m}^{(n)}$ of the number of species of size
$l$ over the enlarged sample of size $n+m$, conditional on the sample
$(X_1,\ldots,X_n)$.

Finally, by combining Theorems~\ref{posteriormarginalnew} and \ref
{posteriormarginaloldnew} with the specific weights \eqref
{eqgnedinweight} it can be easily verified that for any $l\geq1$
\[
N_{l,m}^{(n)}\stackrel{\mathrm{d}}{\longrightarrow} 0,\qquad
M_{l,m}^{(n)}\stackrel{\mathrm{d}}{\longrightarrow} 0
\]
as $m\rightarrow+\infty$. As in the unconditional case, these limits
are not surprising due to the almost sure finiteness of the number of
blocks of a random partition induced by the Gnedin model.

\section{Genomic applications}\label{sec4}

A Bayesian nonparametric model \eqref{eqexchangeabl}, with $\tilde P$
being a Gibbs-type r.p.m.~with $\sigma>0$, is particularly suited for
inferential problems with a large unknown number of species given it
postulates an infinite number of species. These usually occur in
genomic applications, such as the analysis of Expressed Sequence Tags
(EST), Cap Analysis Gene Expression (CAGE) or Serial Analysis of
Gene Expression (SAGE). See, for example, \cite{Val09,lmp07bmc,Dur09}. The typical situation is as follows: a
sample of size $n$ sequenced from a genomic library is available and
one would like to make predictions, over an enlarged sample of size
$n+m$ and conditionally on the observed sample, of certain quantities
of interest. The most obvious quantity is the number of distinct
species to be observed in the enlarged sample, which represents a
measure of the overall genes variety. The resulting Bayesian
nonparametric estimators proposed in \cite{Lij207,Lij308} have
already been integrated into the web server RichEst$^\copyright$ \cite
{Dur09}. However, estimators for the overall genes variety are
certainly useful but necessarily need to be complemented by an
effective analysis of the so-called ``rare genes variety'' (see, e.g.,
\cite{Val09}). Therefore, from an applied perspective it is important
to devise estimators of the number of genes that appear only once, the
so-called \textit{unigenes} or, more generally, of the number of genes
that are observed with frequency less than or equal to a specific
abundance threshold $\tau$. The results deduced in the present paper
perfectly fit these needs. Indeed, conditional on an observed sample of
size $n$, the quantity $\hat{M}_{1,m}^{(n)}=\E[M_{l,m}^{(n)}]$ is a
Bayesian estimator of the number of genes that will appear only once in
a sample of size $n+m$ and can be easily determined from Theorem~\ref
{posteriormarginaloldnew}. In a similar fashion, having fixed a
threshold $\tau$,
%
%
\begin{equation}
\hat{M}_{\tau}^{(n)}=\sum_{l=1}^\tau\hat{M}_{l,m}^{(n)}
\label{eqrareestimat}
\end{equation}
is a Bayesian estimator of the rare genes variety, namely, the number
of species appearing with frequency less than $\tau$ in a sample of
size $n+m$.

Having laid out the framework and described the estimators to be used,
we now test the proposed methodology on some real genomic data. To this
end, we deal with a widely used EST data set obtained by sequencing a
tomato-flower cDNA library (made from 0--3 mm buds of tomato flowers)
from the Institute for Genomic Research Tomato Gene Index with library
identifier T1526 \cite{quack}. The observed sample consists of $n=2586$
ESTs with $j=1825$ unique genes, whose frequencies can be summarized by
\[
m_{i,2586}=1434,  253,  71,  33,  11,  6,  2,  3,  1,  2,
2,  1,  1,  1,  2,  1,  1
\]
with $i\in\{1,2,\ldots,14\}\cup\{16,  23,  27\}$, which means that we
are observing $1434$ genes which appear once, $253$ genes which
appear twice, etc.

As for the specific model \eqref{eqexchangeabl} we adopt, $\tilde P$
is a $\operatorname{PD}(\sigma,\theta)$ process. The reason we rely on such a
specification is two-fold: on the one hand, it yields tractable
estimators that can be exactly evaluated and, on the other, it is a
very flexible model since it encompasses a wide range of partitioning
structures according as to the value of $\sigma$. On the basis of our
choice of the nonparametric prior, we only need to specify the
parameter vector $(\sigma, \theta)$. This is achieved by adopting an
empirical Bayes procedure \cite{Lij207}: we fix $(\sigma,\theta)$ so
as to maximize \eqref{eqeppf} corresponding to the observed sample
$(j,n_1,\ldots,n_j)$, that is,
%
%
\begin{equation}
\label{eqmax}
(\hat\sigma,\hat\theta)=
\arg\max_{(\sigma,\theta)}\frac{\prod_{i=1}^{j-1}(\theta+i\sigma
)}{(\theta+1)_{n-1}}
\prod_{i=1}^j (1-\sigma)_{n_i-1}.
\end{equation}
The quantities we wish to estimate are $N_\tau^{(n)}=\sum_{l=1}^{\tau}
N_{l,m}$ and $O_\tau^{(n)}=\break\sum_{l=1}^{\tau}O_{l,m}$. These quantities
identify the number of distinct genes with abundances not greater than
$\tau$ or, in genomic terminology, with expression levels not greater
than $\tau$ that are present among the ``new'' genes detected in the
additional sample and the ``old'' genes observed in the basic sample,
respectively. The overall number of rare distinct genes is easily
recovered as $M_\tau^{(n)}=N_\tau^{(n)}+O_\tau^{(n)}$. The
corresponding estimators can be deduced from \eqref
{bayesestimatorfav1}, \eqref{bayesestimatorfav2} and \eqref
{pippo}. In the present genomic context one can reasonably identify the
rare genes as those presenting expression levels less than or equal to
$\tau=3,4,5$, which are the thresholds we employ for our analysis.

We first perform a cross-validation study for assessing the
performance of our methodology when used to predict rare genes
abundance. To this end, $10$ subsamples of size $1000$ have been drawn
without replacement from the available $2586$ EST sample. For each of
the subsamples we have generated, the corresponding values of $(\sigma
,\theta)$ have been fixed according to \eqref{eqmax}. Predictions
have, then, been performed for an additional sample of size $m=1586$,
which corresponds to the remaining observed genes. Table~\ref{tab1} below
reports the true and estimated values for the $O_\tau^{(n)}$, $N_\tau
^{(n)}$ and $M_\tau^{(n)}$ and shows the accurate performance of the
proposed estimators. Such a result is a fortiori appreciable if one
considers that predictions are made over an additional sample of size
larger than $1.5$ times the observed sample.

%
\begin{table}
\caption{Cross-validation study based on subsamples
of size $1000$ and prediction on the remaining $m=1586$ data. The
reported estimated
and true
quantities are the number of rare genes (i.e., with expression levels
less than or equal to
$\tau$, for $\tau=3,4,5$) among the ``old'' genes ($O_\tau^{(n)}$), the
``new'' genes ($N_\tau^{(n)}$) and all genes ($M_\tau^{(n)}$)}\label{tab1}
\begin{tabular*}{\textwidth}{@{\extracolsep{\fill}}lcccccccccc@{}}
\hline
& & \multicolumn{3}{c}{$\bolds{\tau=3, n=1000}$} &
\multicolumn{3}{c}{$\bolds{\tau=4, n=1000}$} & \multicolumn{3}{c@{}}{$\bolds{\tau
=5, n=1000}$}
\\[-6pt]
& & \multicolumn{3}{c}{\hrulefill} &
\multicolumn{3}{c}{\hrulefill} & \multicolumn{3}{c@{}}{\hrulefill}
\\
\multicolumn{1}{@{}l}{\textbf{N.}} & & $\bolds{O_\tau^{(n)}}$ & $\bolds{N_\tau^{(n)}}$ &
$\bolds{M_\tau^{(n)}}$ & $\bolds{O_\tau^{(n)}}$
& $\bolds{N_\tau^{(n)}}$ & $\bolds{M_\tau^{(n)}}$ & $\bolds{O_\tau^{(n)}}$ & $\bolds{N_\tau^{(n)}}$ &
$\bolds{M_\tau^{(n)}}$ \\
\hline
\phantom{0}{1} & est. & 750 & 1010 & 1759 & 777 & 1014 & 1791 &
793& 1016& 1809 \\
& true & 767 & \phantom{0}991 & 1758 & 793 & \phantom{0}998 & 1791 &
803& \phantom{0}999& 1802 \\ [2.5pt]
\phantom{0}{2} & est. & 739 & 1006 & 1744 & 765 & 1010 & 1775
& 781& 1011& 1792 \\
& true & 753 & 1005 & 1758 & 785 & 1006 & 1791
& 794& 1008& 1802 \\[2.5pt]
\phantom{0}{3} & est. & 730 & 1003 & 1733 & 755 & 1007 & 1762
& 770& 1008& 1779 \\
& true & 742 & 1016 & 1758 & 772 & 1019 & 1791
& 783& 1019& 1802 \\[2.5pt]
\phantom{0}{4} & est. & 765 & 1043 & 1807 & 789 & 1047 & 1836
& 804& 1048& 1852 \\
& true & 772 & \phantom{0}986 & 1758 & 800 & \phantom{0}991 & 1791
& 811& \phantom{0}991& 1802 \\[2.5pt]
\phantom{0}{5} & est. & 741 & \phantom{0}971 & 1712 & 771 & \phantom
{0}976 & 1748
& 788& \phantom{0}978& 1766 \\
& true & 761 & \phantom{0}997 & 1758 & 788 & 1003 & 1791
& 797& 1005& 1802 \\[2.5pt]
\phantom{0}{6} & est. & 758 & 1027 & 1785 & 784 & 1031 & 1816
& 800& 1033& 1833 \\
& true & 770 & \phantom{0}988 & 1758 & 798 & \phantom{0}993 & 1791
& 809& \phantom{0}993& 1802 \\[2.5pt]
\phantom{0}{7} & est. & 739 & \phantom{0}997 & 1735 & 766 & 1002 & 1768
& 783& 1003& 1786 \\[2.5pt]
& true & 758 & 1000 & 1758 & 787 & 1004 & 1791
& 796& 1006& 1802 \\[2.5pt]
\phantom{0}{8} & est. & 734 & \phantom{0}984 & 1719 & 763 & \phantom
{0}989 & 1752
& 780& \phantom{0}991& 1770 \\
& true & 747 & 1011 & 1758 & 779 & 1012 & 1791
& 790& 1012& 1802 \\[2.5pt]
\phantom{0}{9} & est. & 729 & \phantom{0}969 & 1698 & 759 & \phantom
{0}974 & 1733
& 777& \phantom{0}975& 1752 \\
& true & 747 & 1011 & 1758 & 779 & 1012 & 1791
& 789& 1013& 1802 \\[2.5pt]
{10} & est. & 757 & 1020 & 1777 & 784 & 1025 & 1809
& 800& 1026& 1826 \\
& true & 774 & \phantom{0}984 & 1758 & 799 & \phantom{0}992 & 1791
& 807& \phantom{0}995& 1802 \\
\hline
\end{tabular*}
\end{table}

We now deal with the whole data set and provide estimates of rare genes
abundance after additional sequencing. To this end, we consider, as
possible sizes of the additional sample,\vadjust{\goodbreak} $m \in\{250,500,750,1000\}$.
As for the prior specification of $(\sigma, \theta)$, the maximization
in \eqref{eqmax} leads to $(\hat\sigma,\hat\theta)=(0.612,  741)$.
The resulting estimates of $O_\tau^{(n,j)}$, $N_\tau^{(n,j)}$ and
$M_\tau^{(n,j)}$ are reported in Table~\ref{tab2}.

\begin{table}
\caption{Estimates for an additional sample
corresponding to
$m \in\{250,500,750,1000\}$ given the observed EST data set of size
$n=2586$ with $j=1825$ distinct genes:
estimates for the number of rare genes (i.e., with expression levels
less than or equal to
$\tau$, for $\tau=3,4,5$) among the ``old'' genes ($O_\tau^{(n)}$), the
``new'' genes ($N_\tau^{(n)}$) and all genes ($M_\tau^{(n)}$)}\label{tab2}
\begin{tabular*}{\textwidth}{@{\extracolsep{\fill}}lccccccccc@{}}
\hline
& \multicolumn{3}{c}{$\bolds{\tau=3}$} &
\multicolumn{3}{c}{$\bolds{\tau=4}$} & \multicolumn{3}{c@{}}{$\bolds{\tau=5}$} \\[-6pt]
& \multicolumn{3}{c}{\hrulefill} &
\multicolumn{3}{c}{\hrulefill} & \multicolumn{3}{c@{}}{\hrulefill}
\\
& \multicolumn{3}{c}{$\bolds{n=2586,}$ $\bolds{j=1825}$} &
\multicolumn{3}{c}{$\bolds{n=2586,}$ $\bolds{j=1825}$} & \multicolumn{3}{c}{$\bolds{n=2586,}$
$\bolds{j=1825}$}\\[-6pt]
& \multicolumn{3}{c}{\hrulefill} & \multicolumn{3}{c}{\hrulefill} &
\multicolumn{3}{c@{}}{\hrulefill} \\
{$\bolds{m}$} & $\bolds{\hat O_\tau^{(n)}}$ & $\bolds{\hat N_\tau^{(n)}}$ & $\bolds{\hat M_\tau
^{(n)}}$ &
$\bolds{\hat O_\tau^{(n)}}$ & $\bolds{\hat N_\tau^{(n)}}$ & $\bolds{\hat M_\tau^{(n)}}$ &
$\bolds{\hat
O_\tau^{(n)}}$ & $\bolds{\hat N_\tau^{(n)}}$ & $\bolds{\hat M_\tau^{(n
)}}$ \\
\hline
\phantom{0}250 & 1745 & 138 & 1882 & 1782 & 138 & 1920 & 1798 & 138 &
1935 \\
\phantom{0}500 & 1730 & 272 & 2002 & 1773 & 272 & 2045 & 1793 & 272 &
2064 \\
\phantom{0}750 & 1715 & 402 & 2117 & 1763 & 402 & 2165 & 1787 & 403 &
2189\\
1000 & 1700 & 529 & 2229 & 1753 & 530 & 2283 & 1780 & 530 &
2310 \\
\hline
\end{tabular*}
\end{table}
%


\section{Proofs}\label{sec5}

We start by providing a lemma concerning the marginal frequency counts
of the partition blocks
induced by Gibbs-type random partition. In addition to the notation
introduced in Section~\ref{sec2}, we define
the following shortened set notation:
\[
A_{n,m}(j,\bm{n},s,k):=\bigl\{K_n=j,\bm{N}=\bm
{n},L_m^{(n)}=s,K_m^{(n)}=k\bigr\}
\]
and
\[
A_{n}(j,\bm{n}):=\{K_n=j,\bm{N}=\bm{n}\}
\]
for any $\bm{n}=(n_1,\ldots,n_j)\in\mathcal{D}_{n,j}$. Further
additional notation will be introduced in
the proofs when necessary.

\begin{lem}\label{lemmamarginals}
Let $(X_{n})_{n\geq1}$ be an exchangeable sequence associated to a
Gibbs-type r.p.m.
For any $x\in\{1,\ldots,j\}$, let $\bm{q}^{(x)}=(q_1,\ldots,q_x)$ with
$1\le q_1<\cdots<q_x\le j$ and define the vector of frequency counts
$\bm{S}_{\bm{q}^{(x)}}:=(S_{q_1},\ldots,S_{q_x}).$
Then,
\begin{eqnarray}
\label{marginalsold}
&&\P[\bm{S}_{\bm{q}^{(x)}}=\bm{s}_{\bm{q}^{(x)}} | A_{n,m}(j,\bm
{n},s,k)]\nonumber\\
&&\qquad=\frac{(m-s)!}{(m-s-|\bm{s}_{\bm{q}^{(x)}}|)!}
\prod_{i=1}^{x}\frac{(n_{q_{i}}-\sigma)_{s_{q_{i} } }}{s_{q_i}!}\\
&&\qquad\quad{}\times\frac{(n-|\bm{n}_{\bm{q}^{(x)}}|-(j-x)\sigma)_{m-s-|\bm
{s}_{\bm
{q}^{(x)}}|} }{ (n-j\sigma)_{m-s}}\nonumber
\end{eqnarray}
for any vector $\bm{s}_{\bm{q}^{(x)}}=(s_{q_1},\ldots,s_{q_x})$ of
nonnegative integers such that $|\bm{s}_{\bm{q}^{(x)}}|=\sum
_{i=1}^xs_{q_i}\le m-s$.
Moreover, for any $y\in\{1,\ldots,k\}$, let $\bm
{r}^{(y)}=(r_1,\ldots
,r_y)$ with $1\le r_1<\cdots<r_y\le k$ and define the vector of
frequency counts $
\bm{S}_{\bm{r}^{(y)}}^*:=(S_{j+r_{1}},\ldots,S_{j+r_{y}})$. Then
\begin{eqnarray}\label{marginalsnew}
&&\P[\bm{S}_{\bm{r}^{(y)}}^*=\bm{s}_{\bm{r}^{(y)}} |
A_{n,m}(j,\bm
{n},s,k)]\nonumber\\
&&\qquad=\frac{s!}{(s-|\bm{s}_{\bm{r}^{(y)}}|)!}
\prod_{i=1}^{y}\frac{(1-\sigma)_{s_{j+r_i}-1}}{s_{j+r_i}!}\\
&&\qquad\quad{}\times\frac{(k-y)!}{k!}\sigma^{y}
\frac{\mathscr{C}(s-|\bm{s}_{\bm{r}^{(y)}}|,k-y;\sigma)}
{\mathscr{C}(s,k;\sigma)}\nonumber
\end{eqnarray}
for any vector $\bm{s}_{\bm{r}^{(y)}}=(s_{j+r_1},\ldots,s_{j+r_y})$ of
positive integers such that $|\bm{s}_{\bm{r}^{(y)}}|=\sum
_{i=1}^{y}s_{j+r_i}\le s$. Moreover, the random
variables $\bm{S}_{\bm{q}^{(x)}}$ and $\bm{S}_{\bm{r}^{(y)}}^*$ are
independent, conditionally on $(K_{n},\bm{N},L_{m}^{(n)},K_{m}^{(n)})$.
\end{lem}

\begin{pf}
We start by recalling some useful conditional formulae for Gibbs-type
random partitions recently obtained
in \cite{Lij308}. In particular, from \cite{Lij308}, Corollary~1,
one has the
conditional probability
\begin{eqnarray}
\label{eqjointdistinctsample}
&&\P\bigl[K_{m}^{(n)}=k,L_{m}^{(n)}=s
|  A_{n}(j,\bm{n})\bigr]
\nonumber
\\[-8pt]
\\[-8pt]
\nonumber
&&\qquad=\frac{V_{n+m,j+k}}{V_{n,j}}\pmatrix{m\vspace*{2pt}\cr
s} (n-j\sigma)_{m-s}
\frac{\mathscr{C}(s,k,\sigma)}{\sigma^k}.
\end{eqnarray}
On the other hand, for any vectors of nonnegative integers $\bm
{s}_{\bm
{q}^{(j)}}=(s_1,\ldots,s_j)$ such that $|\bm{s}_{\bm{q}^{(j)}}|=m-s$,
and for any vector of positive integers $\bm{s}_{\bm
{r}^{(k)}}=(s_{j+1},\ldots,s_{j+k})$ such that $|\bm{s}_{\bm
{r}^{(k)}}|=s$, according to \cite{Lij308}, equation (28), the
expression
\begin{equation}
\label{eqjointconditional}
\frac{V_{n+m,j+k}}{V_{n,j}}
\prod_{i=1}^{j}(n_i-\sigma)_{s_{i}} \prod_{\ell=1}^{k}(1-\sigma
)_{s_{j+\ell}-1}
\end{equation}
is the conditional probability, given $A_n(j,\mathbf{n})$, of observing a sample
$X_{n+1},\ldots,\break X_{n+m}$ such
that: (i) $L^{(n)}_{m} = s$ elements generate $K^{(n)}_{m} = k$ new distinct species with frequencies $\mathbf{s}_{\mathbf{r}(k)}$ and
(ii) the remaining $m- s$ elements coincide with any of the~$j$ distinct species in $X_1,\ldots,
X_n$
and display a vector of frequencies $\mathbf{s}_{\mathbf{q}^{(j)}}$. Hence, (\ref{eqjointconditional}) determines the conditional probability
distribution of $(\mathbf{S}_{\mathbf{q}^{(j)}} , \mathbf{S}_{\mathbf{r}^{(k)}}^{\ast}
,L^{(n)}_{m}, K^{(n)}_{m})$, given $X_1,\ldots,X_n$.
A combination of \eqref{eqjointdistinctsample} and \eqref
{eqjointconditional} implies that
\begin{equation}
\label{eqdistribution}
\frac{\sigma^k  \prod_{i=1}^{j}(n_i-\sigma)_{s_{q_i}-1} \prod
_{\ell
=1}^{k}(1-\sigma)_{s_{j+r_\ell}-1}}{{m\choose
s}(n-j\sigma)_{m-s} \mathscr{C}(s,k,\sigma)}
\end{equation}
yields the conditional probability, given $A_{n,m}(j,\mathbf{n}, s, k)$, of observing a sample $X_{n+1},\ldots,X_{n+m}$
such that: (i) the $k$ new distinct species featured by $s$ of the~$m$ observations have frequencies
$\mathbf{s}_{\mathbf{r}(k)}$ and (ii) the remaining $m- s$ elements coincide with any of the $j$ distinct species in
$X_1,\ldots,X_n$ and display a vector of frequencies $\mathbf{s}_{\mathbf{q}^{(j)}}$. Hence, (\ref{eqdistribution}) determines the conditional
probability distribution of $(\mathbf{S}_{\mathbf{q}^{(j)}},\mathbf{S}_{\mathbf{r}^{(k)}}^{\ast} )$, given
$(X_{1},\ldots,X_n,L_{m}^{(n)}
,K^{(n)}_{m})$.
Consider now the set $\mathcal{I}_{j,x}:=\{1,\ldots,j\}\setminus\{
q_{1},\ldots,q_{x}\}$ and the corresponding partition set
defined as follows:
\[
\mathcal{D}_{m-s-s^*,j-x}^{(0)}:=\biggl\{(s_{i},i\in\mathcal
{I}_{j,x}) \dvtx
s_{i}\geq0\mbox{ and }\sum_{i\in\mathcal
{I}_{j,x}}s_{i}=m-s-s^*\biggr\},
\]
where we set $s^*:=\sum_{i=1}^{x}s_{q_i}$. In a similar vein, let us
introduce the set
$\mathcal{I}_{k,y}:=\{1,\ldots,k\}\setminus\{r_{1},\ldots,r_{y}\}$
and the corresponding partition set defined as follows:
\[
\mathcal{D}_{s-s^{**},k-y}:=\biggl\{(s_{j+i},i\in\mathcal
{I}_{k,y}) \dvtx
s_{j+i}>0\mbox{ and }
\sum_{i\in\mathcal{I}_{k,y}}s_{j+i}=s-s^{**}\biggr\},
\]
where we set $s^{**}:=\sum_{i=1}^{y}s_{j+r_{i}}$. By virtue of \cite{Cha05}, equation (2.6.1), one can write
\begin{eqnarray}
\label{eqsommacharalambo}
&&\frac{1}{(k-y)!}\sum_{\mathcal{D}_{s-s^{**},k-y} } s!
\prod_{i=1}^{k}\frac{(1-\sigma)_{s_{j+i}-1}}{s_{j+i}!}
\nonumber
\\[-8pt]
\\[-8pt]
\nonumber
&&\qquad=\frac{s!}{(s-s^{**})!\prod_{i=1}^ys_{r_i}!}
\frac{\mathscr{C}(s-s^{**},k-y,\sigma)
}{\sigma^{k-y}}
\end{eqnarray}
and, by virtue of \cite{Lij308}, Lemma (A.1), one can write
\begin{eqnarray}
\label{eqlpw}
&&\sum_{\mathcal{D}_{m-s-s^*,j-x}^{(0)}}\pmatrix{m-s\vspace*{2pt}\cr s_{1},\ldots
,s_{j} }
\prod_{i=1}^{j}(1-\sigma)_{n_{i}+s_{i}-1 }
\nonumber
\\
&&\qquad=\frac{(m-s)!(n^*-(j-x)\sigma)_{m-s-s^*}}{(m-s-s^*)! \prod_{i=1}^x
s_{q_i}!}\\
&&\quad\qquad{}\times\prod_{i=1}^x(1-\sigma)_{n_{q_i}+s_{q_i}-1}\prod_{\ell\in\mathcal
{I}_{j,x}}
(1-\sigma)_{n_\ell-1},\nonumber
\end{eqnarray}
where we set $n^*:=\sum_{i\in\mathcal{I}_{j,x}}n_i=n-\sum_{i=1}^x n_{q_i}$.
A simple application of the identities \eqref{eqsommacharalambo} and
\eqref{eqlpw} to the conditional probability \eqref{eqdistribution}
proves both the conditional independence between $\bm{S}_{\bm
{q}^{(x)}}$ and $\bm{S}_{\bm{r}^{(y)}}^*$ and the two expressions in
\eqref{marginalsold} and \eqref{marginalsnew}.
\end{pf}

\subsection{\texorpdfstring{Proof of Proposition~\protect\ref{priormarginal}}
{Proof of Proposition 1}}\label{sec5.1} For any $n\geq1$
and $1\leq j\leq n$ let $\mathcal{M}_{n,j}$ be the partition set of
$\N
_n$ containing all the vectors $\bm{m}_n=(m_{1},\ldots,m_{n})\in\{
0,1,\ldots,n\}^n$ such that\vadjust{\goodbreak} $\sum_{i=1}^{n}m_i=j$ and $\sum
_{i=1}^{n}im_{i}=n$. Hence, resorting to the probability distribution
\eqref{eqgibbssampling}, one obtains for any $r\geq1$
\begin{eqnarray*}
\E\bigl[(M_{l,n})_{[r]}\bigr]
&=&n!\sum_{j=1}^{n}V_{n,j}\sum_{\bm{m}_n\in\mathcal
{M}_{n,j}}(m_{l})_{[r]}\prod_{i=1}^{n}\biggl(\frac{(1-\sigma
)_{i-1}}{i!}\biggr)^{m_{i}}\frac{1}{m_{i}!}\\
&=&n!\sum_{j=1}^{n}V_{n,j}\sum_{\bm{m}_{n}\in\mathcal
{M}_{n,j}}
\biggl(\frac{(1-\sigma)_{l-1}}{l!}\biggr)^{m_{l}}\frac{1}{(m_{l}-r)!}
\\
& &{} \times
 \prod_{1\leq i\neq l\leq n }\biggl(\frac{(1-\sigma
)_{i-1}}{i!}
\biggr)^{m_{i}}\frac{1}{m_{i}!}\\
&=&n!\biggl(\frac{(1-\sigma)_{l-1}}{l!}\biggr)^{r} \sum
_{j=1}^{n}V_{n,j}
 \sum_{\bm{m}_{n-rl}\in\mathcal{M}_{n-rl,j-r}}
\prod_{i=1}^{n-rl}\biggl(\frac{(1-\sigma)_{i-1}}{i!}
\biggr)^{m_{i}}\hspace*{-1pt}\frac
{1}{m_{i}!}.
\end{eqnarray*}
Finally, a direct application of \cite{Cha05}, equation (2.82),
implies the following identity:
\[
\sum_{\bm{m}_{n}\in\mathcal{M}_{n-rl,j-r}}
\prod_{i=1}^{n}\biggl(\frac{(1-\sigma)_{i-1}}{i!}
\biggr)^{m_{i}}\frac{1}{m_{i}!}
=\frac{(n)_{[lr]}}{n!\sigma^{j-r}}\mathscr{C}(n-lr,j-r;\sigma),
\]
and the proof is completed.

\subsection{\texorpdfstring{Proof of Theorem~\protect\ref{posteriormarginalold}}
{Proof of Theorem 1}}\label{sec5.2} According to
the definition of the random variable $O_{l,m}$ in \eqref
{eqoldnewsmallblocks},
for any $r\geq1$ one can write
\begin{eqnarray*}
\E\bigl[\bigl(O_{l,m}^{(n)}\bigr)^{r}\bigr]&=&
\sum_{s=0}^{m}\sum_{k=0}^{s}\P\bigl[L_{m}^{(n)}=s,K_{m}^{(n)}=k |
A_n(j,\bm{n})\bigr]
\\
&&{}\times\E\Biggl[\Biggl(\sum_{i=1}^j\indic_{l}(n_i+S_i)\Biggr)^{r}
\Big|  A_{n,m}(j,\bm{n},s,k)\Biggr].
\end{eqnarray*}
It can be easily verified that a repeated application of the binomial
expansion implies the following identity:
\begin{eqnarray}\label{eqidentity1}
&&\biggl(\sum_{i=1}^j\indic_{\{l\}}(n_i+S_i)\biggr)^{r}\nonumber\\
&&\qquad=
\sum_{x=1}^{j}\sum_{i_{1}=1}^{r-1}\sum_{i_{2}=1}^{i_{2}-1}\cdots
\sum_{i_{x-1}=1}
^{i_{x-2}-1}
\pmatrix{r\vspace*{2pt}\cr i_{1}}\pmatrix{i_{1}\vspace*{2pt}\cr i_{2}}\cdots\pmatrix{i_{x-2}\vspace*{2pt}\cr i_{x-1}}
\\
&&\quad\qquad\hspace*{95pt}{}\times \sum_{\bm{c}^{(x)}\in
\mathcal{C}_{j,x}}\prod_{t=1}^{x}\bigl(\indic_{\{l\}}(n_{c_{t}}+S_{c_{t}})
\bigr)^{i_{x-t}-i_{x-t+1}}\nonumber
\end{eqnarray}
provided $i_0\equiv r$. Observe that the previous sum can be expressed
in terms of Stirling numbers of the second kind $S(n,m)$; indeed, since
$m!S(n,m)$ is the number of ways of distributing $n$ distinguishable
objects into $m$ distinguishable groups, one has
%
%
\begin{equation}
\frac{1}{m!}\sum_{i_{1}=1}^{n-1}\sum_{i_{2}=1}^{i_{1}-1}\cdots\sum
_{i_{m-1}=1}^
{i_{m-2}-1}\pmatrix{n\vspace*{2pt}\cr i_{1} }\pmatrix{i_{1}\vspace*{2pt}\cr i_{2}
}\cdots\pmatrix{i_{m-2}\vspace*{2pt}\cr
i_{m-1}}=S(n,m),
\label{eqstirlingbinom}
\end{equation}
for any $n\ge1$ and $1\le m\le n$. In particular, combining the
identity \eqref{eqidentity1} with \eqref{eqstirlingbinom}, one obtains
\begin{eqnarray}
\label{eqidentity1bis}
&&\E\bigl[\bigl(O_{l,m}^{(n)}\bigr)^{r}  \big|
L_{m}^{(n)}=s,K_{m}^{(n)}=k\bigr]
\nonumber
\\[-8pt]
\\[-8pt]
\nonumber
&&\qquad=\sum_{x=1}^{j\wedge r}S(r,x)x! \sum_{\bm{c}^{(x)}\in\mathcal{C}_{j,x}}\P[\bm
{S}_{\bm{c}^{(x)}}=
l\bm{1}_{x}-\bm{n}_{\bm{c}^{(x)}} | A_{n,m}(j,\bm{n},s,k)],
\end{eqnarray}
where we set $\mathbf{1}_{x}:=(1,\ldots,1)$ and $\bm{n}_{\bm
{c}^{(x)}}=(n_{c_1},\ldots,n_{c_x})$. In \eqref{eqidentity1bis} the
bound $j\wedge r$ on the sum over the index $x$ is motivated by the
fact that $S(r,x)=0$ if $x>r$. Hence, the identity \eqref{eqidentity1bis}
combined with \eqref{marginalsold} yields the following expression:
\begin{eqnarray}\label{eqformula11}
&&\E\bigl[\bigl(O_{l,m}^{(n)}\bigr)^{r}  |
L_{m}^{(n)}=s,K_{m}^{(n)}=k\bigr]\nonumber\\
&&\qquad=\sum_{x=1}^{j\wedge r}S(r,x)x! \sum_{\bm{c}^{(x)}\in\mathcal{C}_{j,x}}
\frac{(m-s)!}{(m-s-xl+|\bm{n}_{\bm{c}^{(x)}}|)!} \prod_{i=1}^x
\frac{(n_{c_i}-\sigma)_{l-n_{c_i}}}{(l-n_{c_i})!}\\
&&\hspace*{122pt}{}\times
\frac{(n-|\bm{n}_{\bm{c}^{(x)}}|-(j-x)\sigma
)_{m-s-xl+|\bm
{n}_{\bm{c}^{(x)}}|}}
{(n-j\sigma)_{m-s}}.\nonumber
\end{eqnarray}
Observe that in \eqref{eqformula11} the sum over the index $x$, for
$x=1,\ldots,j\wedge r$, is equivalent to a sum over the index $x$ for
$x=1,\ldots,r$. Indeed, if $j>r$, then the sum over the index $x$ is
nonnull for $x=1,\ldots,r$ because $S(r,x)=0$ for any $x=r+1,\ldots
,j$; on the other hand, if $j<r$, then the sum over the index $x$ is
nonnull for $x=1,\ldots,j$ because the set $\mathcal{C}_{j,x}$ is
empty for any $x=j+1,\ldots,r$. Accordingly, resorting to~\cite{Lij308}, Corollary~1, one can rewrite the expected value above
as
\begin{eqnarray*}
\E\bigl[\bigl(O_{l,m}^{(n)}\bigr)^{r}\bigr]&=&
\sum_{s=0}^{m}\sum_{k=0}^{s}\frac{V_{n+m,j+k}}{V_{n,j}}\pmatrix{m\vspace*{2pt}\cr
s}\frac
{\mathscr{C}(s,k;\sigma)}{\sigma^{k}}\sum_{x=1}^{r}S(r,x)x!\\
&&{}\times\sum_{\bm{c}^{(x)}\in\mathcal{C}_{j,x}}\frac
{(m-s)!}{(m-s-xl+|\bm{n}_{\bm{c}^{(x)}}|)!}
 \prod_{i=1}^{x}\frac{(n_{c_{i}}-\sigma
)_{l-n_{c_{i}}}}{(l-n_{c_i})!}\\
&&{}\times
 \bigl(n-|\bm{n}_{\bm{c}^{(x)}}|-(j-x)\sigma\bigr)_{m-s-xl+|\bm
{n}_{\bm{c}^{(x)}}|}\\
&=&
\sum_{x=1}^{r}S(r,x)x!\sum_{\bm{c}^{(x)}\in\mathcal{C}_{j,x}}\frac
{m!}{(m-xl+|\bm{n}_{\bm{c}^{(x)}}|)!}
 \prod_{i=1}^{x}\frac{(n_{c_{i}}-\sigma
)_{l-n_{c_{i}}}}{(l-n_{c_i})!}\\
&&{}\times\sum_{k=0}^{m-xl+|\bm{n}_{\bm{c}^{(x)}}|}\frac
{V_{n+m,j+k}}{V_{n,j}} \sigma^{-k}
\sum_{s=k}^{m-xl+|\bm{n}_{\bm{c}^{(x)}}|}
\pmatrix{m-xl+|\bm{n}_{\bm{c}^{(x)}}|\vspace*{2pt}\cr s}\\
&&{}\times\bigl(n-|\bm{n}_{\bm
{c}^{(x)}}|-(j-x)\sigma
\bigr)_{m-xl+|\bm{n}_{\bm{c}^{(x)}}|-s}
\mathscr{C}(s,k;\sigma)\\
&=&\sum_{x=1}^{r}S(r,x)x!\sum_{\bm{c}^{(x)}\in\mathcal
{C}_{j,x}}\frac
{m!}{(m-xl+|\bm{n}_{\bm{c}^{(x)}}|)!}
 \prod_{i=1}^{x}\frac{(n_{c_{i}}-\sigma
)_{l-n_{c_{i}}}}{(l-n_{c_i})!}\\
&&{}\times\sum_{k=0}^{m-xl+|\bm{n}_{\bm{c}^{(x)}}|}\frac
{V_{n+m,j+k}}{V_{n,j}}
\\
&&{}\times\frac{\mathscr{C}(m-xl+|\bm{n}_{\bm{c}^{(x)}}|,k;\sigma
,-n+|\bm
{n}_{\bm{c}^{(x)}}|
+(j-x)\sigma)}{\sigma^k},
\end{eqnarray*}
where the last equality follows from \cite{Cha05}, equation (2.56).
The proof of \eqref{eqmarginalposteriorold2} is, thus, completed by
using the relation between the $r$th moment with the $r$th factorial
moment.

\subsection{\texorpdfstring{Proof of Theorem~\protect\ref{posteriormarginalnew}}
{Proof of Theorem 2}}\label{sec5.3} The proof is
along lines similar to the proof of Theorem~\ref{posteriormarginalold}.
In particular, it can be easily verified that a repeated application of
the binomial expansion implies the following identity:
\begin{eqnarray*}
\Biggl(\sum_{i=1}^k\indic_{\{l\}}(S_{j+i})\Biggr)^{r}&=&
\sum_{y=1}^{k}\sum_{i_{1}=1}^{r-1}\sum_{i_{2}=1}^{i_{2}-1}\cdots
\sum
_{i_{y-1}=1}^{i_{y-2}-1}
\pmatrix{r\vspace*{2pt}\cr i_{1}}\pmatrix{i_{1}\vspace*{2pt}\cr i_{2}}\cdots\pmatrix{i_{y-2}\vspace*{2pt}\cr i_{y-1}}\\
&&\hspace*{95pt}{}\times \sum_{\bm{c}^{(y)}\in
\mathcal{C}_{k,y}}\prod_{t=1}^{y}\bigl(\indic_{\{l\}}(S_{j+c_{t}})
\bigr)^{i_{y-t}-i_{y-t+1}}.
\end{eqnarray*}
Hence, according to the definition of the random variable $N_{l,m}$ in
\eqref{eqoldnewsmallblocks} and by combining the identity \eqref
{eqstirlingbinom} with \eqref{marginalsnew}, one has
\begin{eqnarray}
\label{identity2}
&&\E
 \bigl[\bigl(N_{l,m}^{(n)}\bigr)^{r} |
L_{m}^{(n)}=s,K_{m}^{(n)}=k\bigr]\nonumber\\
&&\qquad=\sum_{y=1}^{k}S(r,y)y!\pmatrix{k\vspace*{2pt}\cr y}\P[\bm{S}_{\bm
{c}^{(y)}}^*=l\mathbf{1}_{y} | A_{n,m}(j,\bm{n},s,k)]\\
&&\qquad=\sum_{y=1}^{k}S(r,y)\frac{s!}{(s-yl)!}  \frac{[\sigma
(1-\sigma
)_{l-1}]^y}{(l!)^y}
\frac{\mathscr{C}(s-yl,k-y;\sigma)}{\mathscr{C}(s,k;\sigma)},\nonumber
\end{eqnarray}
where we set $\bm{1}_y:=(1,\ldots,1)$. Hence, \eqref{identity2}
combined with \eqref{eqjointdistinctsample} leads to the
following expression:
\begin{eqnarray}\label{identity21}
&&\E\bigl[\bigl(N_{l,m}^{(n)}\bigr)^{r}\bigr]\nonumber\\
&&\qquad=\sum_{s=0}^{m}\sum_{k=0}^{s}\frac{V_{n+m,j+k}}{V_{n,j}}\pmatrix{m\vspace*{2pt}\cr
s}(n-j\sigma)_{m-s}\\
&&\hspace*{64pt}{}\times \sum_{y=1}^{r\wedge k}S(r,y) \frac{s!}{(s-yl)!}
\frac{[\sigma(1-\sigma
)_{l-1}
]^y}{(l!)^y}
\frac{\mathscr{C}(s-yl,k-y;\sigma)}{\sigma^{k}}.\nonumber
\end{eqnarray}
In \eqref{identity21} note that the sum over the index $y$, for
$y=1,\ldots,k$, is equivalent to a sum over the index $y$ for
$y=1,\ldots,r$. Indeed, if $k>r$, then the sum over the index $y$ is
nonnull for $y=1,\ldots,r$ because $S(r,y)=0$ for any $y=r+1,\ldots
,k$; on the other hand, if $k<r$, then the sum over the index $y$ in
nonnull for $y=1,\ldots,k$ because $\mathscr{C}(s-yl,k-y;\sigma)=0$
for any $y=k+1,\ldots,r$. Based on this, one can rewrite the expected
value above as
\begin{eqnarray*}
\E\bigl[\bigl(N_{l,m}^{(n)}\bigr)^{r}\bigr]
&=&\sum_{y=1}^{r}S(r,y)\frac{[(1-\sigma)_{l-1}]^y}{(l!)^y}
 \sum_{s=yl}^m\pmatrix{m\vspace*{2pt}\cr s}(n-j\sigma)_{m-s}  \frac{s!}{(s-yl)!}\\
&&{} \times \sum_{k=y}^{s} \frac{V_{n+m,j+k}}{V_{n,j}}
\frac{\mathscr{C}(s-yl,k-y;\sigma)}{\sigma^{k-y}}\\
&=&\sum_{y=1}^{r}S(r,y) \frac{[(1-\sigma)_{l-1}
]^y}{(l!)^y}
\sum_{s=0}^{m-yl}
\pmatrix{{m}\vspace*{2pt}\cr{s+yl}}(n-j\sigma)_{m-s-yl}\frac{(s+yl)!}{(s)!}\\
&&{} \times\sum_{k=0}^{s+yl-y}
\sigma^{-k} \frac{V_{n+m,j+k+y}}{V_{n,j}}\mathscr{C}(s,k;\sigma)\\
&=&\sum_{y=1}^{r}S(r,y) \frac{[(1-\sigma)_{l-1}
]^y}{(l!)^y}
\frac{m!}{(m-yl)!}
 \sum_{k=0}^{m-yl}
\sigma^{-k} \frac{V_{n+m,j+k+y}}{V_{n,j}}\\
&&{} \times\sum_{s=k}^{m-yl}\pmatrix{{m-yl}\vspace*{2pt}\cr{s}}(n-j\sigma)_{m-yl-s}
\mathscr{C}(s,k;\sigma)\\
&=&\sum_{y=1}^{r}S(r,y) [(1-\sigma)_{l-1}]^y \frac
{m!}{(l!)^y (m-yl)!}\\
&&{} \times
 \sum_{k=0}^{m-yl}
\sigma^{-k} \frac{V_{n+m,j+k+y}}{V_{n,j}}
{\mathscr{C}(m-yl,k;\sigma,-n+j\sigma)}.
\end{eqnarray*}
The proof of \eqref{eqmarginalposteriornew2} is, thus, completed by
using the relation between the $r$th moment with the $r$th factorial
moment.

\subsection{\texorpdfstring{Proof of Theorem~\protect\ref{posteriormarginaloldnew}}
{Proof of Theorem 3}}\label{sec5.4} The proof
follows from conditional independence between the random variables $\bm
{S}_{\bm{q}^{(x)}}$ and $\bm{S}_{\bm{r}^{(y)}}$, given $(K_n,\bm
{N}_n,L_m^{(n)},K_m^{(n)})$, as stated in Theorem~\ref{priormarginal}.
Indeed, according to the definition of the random variable $M_{l,m}$,
for any $r\geq1$ one can write
\begin{eqnarray}\label{eqm1}
&&\E\bigl[\bigl(M_{l,m}^{(n)}\bigr)^r\bigr]
\nonumber
\\[-8pt]
\\[-8pt]
\nonumber
&&\qquad=\sum_{t=0}^r\pmatrix{{r}\vspace*{2pt}\cr{t}}
\sum_{s=0}^m\sum_{k=0}^s \alpha_t(l) \beta_{r-t}(l) \P
\bigl[L_m^{(n)}=s,  K_m^{(n)}=k | A_n(j,\bm{n})\bigr],
\end{eqnarray}
where
\begin{eqnarray*}
\alpha_t(l)
&:=&\E\bigl[\bigl(O_{l,m}^{(n)}\bigr)^t | L_m^{(n)}=s,
K_m^{(n)}=k\bigr]\\
&=&\sum_{x=1}^{j\wedge t}x!S(t,x)
\sum_{\bm{c}^{(x)}\in\mathcal{C}_{j,x}}\P[
\bm{S}_{\bm{c}^{(x)}}=l\bm{1}_x-\bm{n}_{\bm{c}^{(x)}}
| A_{n,m}(j,\bm{n},s,k)],
\end{eqnarray*}
and
\begin{eqnarray*}
\beta_{r-t}(l)
&:=&\E\bigl[\bigl(N_{l,m}^{(n)}\bigr)^{r-t} | L_m^{(n)}=s,
K_m^{(n)}=k\bigr]\\
&=&\sum_{y=1}^{k\wedge(r-t)}y!S(r-t,y)
\sum_{\bm{c}^{(y)}\in\mathcal{C}_{k,y}}\P[
\bm{S}_{\bm{c}^{(y)}}^*=l\bm{1}_y
| A_{n,m}(j,\bm{n},s,k)].
\end{eqnarray*}
In particular, by combining \eqref{eqm1} with \eqref{eqformula11}
and \eqref{identity2}, one has
\begin{eqnarray*}
&&\E\bigl[\bigl(M_{l,m}^{(n)}\bigr)^r\bigr]\\
&&\qquad=\sum_{t=0}^r\pmatrix{{r}\vspace*{2pt}\cr{t}}\sum_{s=0}^{m}\sum_{k=0}^{s}\P
\bigl[L_m^{(n)}=s,  K_m^{(n)}=k | A_n(j,\bm{n})\bigr]\\
&&\qquad\quad{} \times\sum_{x=1}^{t}S(t,x)x!\sum_{\bm{c}^{(x)}\in
\mathcal
{C}_{j,x}}\frac{(m-s)!}{(m-s-xl+|\bm{n}_{\bm{c}^{(x)}}|)!} \prod
_{i=1}^{x}\frac{(n_{c_i}-\sigma)_{l-n_{c_i}}}{(l-n_{c_i})!}\\
&&\qquad\quad{}  \times\frac{(n-|\bm{n}_{\bm
{c}^{(x)}}|-(j-x)\sigma
)_{m-s-xl+|\bm{n}_{\bm{c}^{(x)}}|}}{(n-j\sigma)_{m-s}}\\
& &\qquad\quad{}  \times\sum_{y=1}^{r-t}S(r-t,y) \frac
{s!}{(s-yl)!}\frac{[\sigma(1-\sigma)_{l-1}]^y}{(l!)^y}\\
&&\hspace*{28pt}\qquad\quad{}\times\frac{\mathscr{C}(s-yl,k-y;\sigma)}{\mathscr{C}(s,k;\sigma)}.
\end{eqnarray*}
Using the same arguments applied in the last part of Theorems~\ref
{posteriormarginalnew} and~\ref{posteriormarginalold}, the
expression \eqref{eqm1} combined with \eqref
{eqjointdistinctsample} leads to the following:
\begin{eqnarray}\label{eqm3}
\E\bigl[\bigl(M_{l,m}^{(n)}\bigr)^r\bigr]
&=&\sum_{t=0}^{r}\pmatrix{r\vspace*{2pt}\cr t}\sum_{x=1}^{t}S(t,x)\sum
_{y=1}^{r-t}S(r-t,y)x!\frac{[(1-\sigma)_{l-1}]^{y}}{(l!)^{y}}\nonumber\\[-2pt]
&&{}\times\sum_{\bm{c}^{(x)}\in\mathcal{C}_{j,x}}\frac
{m!}{(m-xl-yl+|\bm{n}_{\bm{c}^{(x)}}|)!}\prod_{i=1}^{x}\frac
{(n_{c_{i}}-\sigma)_{l-n_{c_{i}}}}{(l-n_{c_{i}})!}\nonumber\\[-9pt]\\[-9pt]
&&{}\times\sum_{k=0}^{m-xl-yl+|\bm{n}_{\bm
{c}^{(x)}}|}\frac
{V_{n+m,j+k+y}}{V_{n,j}}\nonumber\\[-2pt]
&&{}\times\frac{\mathscr{C}(m-xl-yl+|\bm
{n}_{\bm
{c}^{(x)}}|,k;\sigma,-n+|\bm{n}_{\bm{c}^{(x)}}|+(j-x)\sigma
)}{\sigma^{k}}.\nonumber
\end{eqnarray}
The expression in \eqref{eqm3} can be further simplified by applying
well-known properties of the Stirling numbers
of the second kind. In particular, according to the identity
\[
S(r,y+x)\pmatrix{y+x\vspace*{2pt}\cr x}=\sum_{t=x}^{r-y}\pmatrix{r\vspace*{2pt}\cr t}S(t,x)S(r-t,y)
\]
(see \cite{Cha05}, Chapter 2), one can write
\begin{eqnarray*}
\E\bigl[\bigl(M_{l,m}^{(n)}\bigr)^r\bigr]&=&\sum_{x=0}^{r}\sum_{y=0}^{r-x}S(r,y+x)\pmatrix{y+x\vspace*{2pt}\cr x}x!\frac
{[(1-\sigma)_{l-1}]^{y}}{(l!)^{y}}\\[-2pt]
& &{}\times\sum_{\bm{c}^{(x)}\in\mathcal{C}_{j,x}}\frac
{m!}{(m-xl-yl+|\bm{n}_{\bm{c}^{(x)}}|)!}\prod_{i=1}^{x}\frac
{(n_{c_{i}}-\sigma)_{l-n_{c_{i}}}}{(l-n_{c_{i}})!}\\[-2pt]
&&{}  \times\sum_{k=0}^{m-xl-yl+|\bm{n}_{\bm
{c}^{(x)}}|}\frac{V_{n+m,j+k+y}}{V_{n,j}}\\[-2pt]
&&{}   \times\frac{\mathscr{C}(m-xl-yl+|\bm
{n}_{\bm
{c}^{(x)}}|,k;\sigma,-n+|\bm{n}_{\bm{c}^{(x)}}|+(j-x)\sigma
)}{\sigma
^{k}}\\[-2pt]
&=&\sum_{y=0}^{r}S(r,y) \sum_{x=0}^{y}\pmatrix{y\vspace*{2pt}\cr x}x!\frac{[(1-\sigma
)_{l-1}]^{y-x}}{(l!)^{y-x}}\\[-2pt]
&&{} \times\sum_{\bm{c}^{(x)}\in\mathcal{C}_{j,x}}\frac
{m!}{(m-yl+|\bm
{n}_{\bm{c}^{(x)}}|)!}\prod_{i=1}^{x}\frac{(n_{c_{i}}-\sigma
)_{l-n_{c_{i}}}}{(l-n_{c_{i}})!}\\[-2pt]
&&{}  \times\sum_{k=0}^{m-yl+|\bm{n}_{\bm{c}^{(x)}}|}\frac
{V_{n+m,j+k+y-x}}{V_{n,j}}\\[-2pt]
&&{}   \times\frac{\mathscr{C}(m-yl+|\bm{n}_{\bm
{c}^{(x)}}|,k;\sigma,-n+|\bm{n}_{\bm{c}^{(x)}}|+(j-x)\sigma
)}{\sigma^{k}}.
\end{eqnarray*}
The proof of \eqref{eqmarginalposterioroldnew2} is thus completed
by using the relation between the $r$th moment with the $r$th factorial
moment.

\subsection{Proofs for the Dirichlet process}\label{sec5.5}

\subsubsection{\texorpdfstring{Proof of Propositions~\protect\ref{prpdirprior} and \protect\ref{prpdir2}}
{Proof of Propositions 2 and 3}}\label{sec5.5.1}
The distribution of $M_{l,n}$ is determined by its factorial moments as
\begin{eqnarray*}
\P[M_{l,n}=m_{l}]
&=&\frac{\indic_{\{1,\ldots,n\}}(m_ll)}{m_l!} \sum_{k=m_l}^{n}\frac
{(-1)^{k-m_l}}{(k-m_l)!}
\E\bigl[(M_{l,n})_{[k]}\bigr]\\
&=&\frac{n!}{m_l!(\theta)_n} \sum_{k=m_l}^{[n/l]}\frac
{(-1)^{k-m_l}}{(k-m_l)!}
\frac{(\theta)_{[n-kl]}}{l^k(n-kl)!}
\end{eqnarray*}
and, from this, \eqref{eqmarginalpd} easily follows. On the other
hand, Proposition~\ref{prpdir2} is a trivial consequence of \eqref
{eqestimold} and \eqref{eqestimnew}.


\subsection{Proofs for the Pitman model}\label{sec5.6}

\subsubsection{\texorpdfstring{Proof of Proposition~\protect\ref{prp2pd0}}
{Proof of Proposition 4}}\label{sec5.6.1} This again follows
from the application of the sieve formula, as discussed in the proof of
Proposition~~\ref{prpdirprior}.

\subsubsection{\texorpdfstring{Proof of Proposition~\protect\ref{prppdold}}
{Proof of Proposition 5}}\label{sec5.6.2} From Theorem~\ref
{posteriormarginalold} one finds that
\begin{eqnarray*}
\E\bigl[\bigl(O_{l,m}^{(n)}\bigr)_{[r]}\bigr]&=&\frac
{r!m!}{(\theta
+n)_m}
\sum_{\bm{c}^{(r)}\in\mathcal{C}_{j,r}}\frac{1}{(m-rl+|\bm
{n}_{\bm{c}^{(r)}}|)!}
 \prod_{i=1}^r\frac{(n_{c_i}-\sigma)_{l-n_{c_i}}}{(l-n_{c_i})!}\\
&&{}\times \sum_{k=0}^{m-rl+|\bm{n}_{\bm{c}^{(r)}}|}\biggl(\frac
{\theta
}{\sigma}+j\biggr)_k
\\
&&{}\times\Ccr\bigl(m-rl+|\bm{n}_{\bm{c}^{(r)}}|,k;\sigma,-n+|\bm{n}_{\bm
{c}^{(r)}}|+(j-r)\sigma\bigr).
\end{eqnarray*}
By definition,
\[
\sum_{k=0}^n\Ccr(n,k;\sigma,\gamma) (t)_k=(\sigma t-\gamma)_n
\]
and this entails
\begin{eqnarray*}
\E\bigl[\bigl(O_{l,m}^{(n)}\bigr)_{[r]}\bigr]&=&\frac
{r!m!}{(\theta
+n)_m}
\sum_{\bm{c}^{(r)}\in\mathcal{C}_{j,r}}\frac{1}{(m-rl+|\bm
{n}_{\bm{c}^{(r)}}|)!}
 \prod_{i=1}^r\frac{(n_{c_i}-\sigma)_{l-n_{c_i}}}{(l-n_{c_i})!}\\
&&\hspace*{74pt}{}\times (\theta+n-|\bm{n}_{\bm{c}^{(r)}}|+r\sigma)_{m-rl+|\bm
{n}_{\bm
{c}^{(r)}}|}.
\end{eqnarray*}
The usual application of the sieve formula yields \eqref{eqold2pd}.

\subsubsection{\texorpdfstring{Proof of Proposition~\protect\ref{prp2pd2}}
{Proof of Proposition 6}}\label{sec5.6.3} Follows from
Theorem~\ref{posteriormarginalnew}, along the same lines as in the
proof of Proposition~3.4.


\subsubsection{\texorpdfstring{Proof of Theorem~\protect\ref{limitbayesestimatoruniqueness}}
{Proof of Theorem 4}}\label{sec5.6.4}
Our strategy will consist in examining the asymptotic behavior of the
$r$th moments of $N_{l,m}^{(n)}$ and of $M_{l,m}^{(n)}$, for any $r\ge
1$, as $m$ increases. To this end, it is worth referring to the
following decomposition that implicitly follows from the proof of
Theorem~\ref{posteriormarginaloldnew}. Indeed, it can be seen that
\[
\E\bigl[\bigl(M_{l,m}^{(n)}\bigr)^{r}\bigr]=\E\bigl[\bigl(O_{l,m}^{(n)}\bigr)^{r}\bigr]+\E
\bigl[\bigl(N_{l,m}^{(n)}\bigr)^{r}\bigr]
+\sum_{i=1}^{r-1}\pmatrix{r\vspace*{2pt}\cr
i}\mathcal{B}^{(i)}(\sigma,n,j,\bm{n},m),
\]
where
\begin{eqnarray*}
\E\bigl[\bigl(O_{l,m}^{(n)}\bigr)^{r}\bigr]
&=&\frac{m!}{(\theta+n)_m}
\sum_{x=1}^{j\wedge r}x!S(r,x)\sum_{\bm{c}^{(x)}\in\mathcal{C}_{j,x}}
\prod_{r=1}^{x}\frac{(n_{c_r}-\sigma)_{l-n_{c_r}}}{(l-n_{c_r})!}\\
&&{} \times
\frac{(\theta+n-|\bm{n}_{\bm{c}^{(x)}}|+x\sigma)_{m-xl+|\bm
{n}_{\bm{c}^{(x)}}|}}
{
(m-xl+|\bm{n}_{\bm{c}^{(x)}}|)!},\\
\E\bigl[\bigl(N_{l,m}^{(n)}\bigr)^{r}\bigr]
&=&
\frac{m!}{(\theta+n)_m}
\sum_{y=1}^{[m/l]\wedge r}S(r,y) \frac{\sigma^{y}[(1-\sigma
)_{l-1}]^y}{(l!)^y}
 \biggl(j+\frac{\theta}{\sigma}\biggr)_{y }\\
& &{}\times
\frac{(\theta+n+y\sigma)_{m-ly }}{(m-yl)!},\\
\mathcal{B}^{(i)}(\sigma,n,j,\bm{n},m)
&=&\frac{m!}{(\theta+n)_m} \sum_{x=1}^{j\wedge i}x!S(i,x)
\sum_{\bm{c}^{(x)}\in\mathcal{C}_{j,x}}\prod_{r=1}^{x}
\frac{(n_{c_{r}}-\sigma)_{l-n_{c_{r}}}}{(l-n_{c_r})!}\\
&&{}  \times
\sum_{y=1}^{m\wedge
(r-i)}S(r-i,y)\frac{\sigma^{y}[(1-\sigma)_{l-1}]^y}{(l!)^y}
\biggl(j+\frac{\theta}{\sigma}\biggr)_{y}\\
& &{}    \times
 \frac{(\theta+n-|\bm{n}_{\bm{c}^{(x)}}|+\sigma
x)_{m-yl-xl+n_{c_{i}}+|\bm{n}_{\bm{c}^{(x)}}|
}}{(m-yl-xl+|\bm{n}_{\bm{c}^{(x)}}|)!}.
\end{eqnarray*}
By virtue of Stirling's approximation
formula, one has, as $m\rightarrow+\infty$,
\begin{eqnarray*}
&&\E\bigl[\bigl(O_{l,m}^{(n)}\bigr)^{r}\bigr] \\
&&\qquad\sim
m^{-\theta-n+1} \Gamma(\theta+n)\\
&&\qquad\quad{}\times\sum_{x=1}^{j\wedge r} \sum_{\bm{c}^{(x)}\in\mathcal{C}_{j,x}}
x! S(r,x)
 \frac{m^{\theta+n-|\bm{n}_{\bm{c}^{(x)}}|-1+x\sigma
}}{\Gamma
(\theta+n-
|\bm{n}_{\bm{c}^{(x)}}|+x\sigma}
\prod_{t=1}^x\frac{(n_{c_t}-\sigma)_{l-n_{c_t}}}{(l-n_{c_t})!},
\end{eqnarray*}
where $a_n\sim b_n$ means that $a_n/b_n\to1$, as $n\to\infty$. The
term that asymptotically dominates the right-hand side of the
asymptotic equivalence above, as $m\to\infty$, can be bounded by
\[
m^{(j\wedge r)\sigma-|\bm{n}_{\bm{c}^{(j\wedge r)}}|}
\frac{\Gamma(\theta+n) (j\wedge r)!S(r,(j \wedge r))}{\Gamma
(\theta+n-
|\bm{n}_{\bm{c}^{(j\wedge r)}}|+(j\wedge r)\sigma)}\prod
_{t=1}^{(j\wedge r)}
\frac{(n_{c_t}-\sigma)_{l-n_{c_t}}}{(l-n_{c_t})!}.
\]
Since $|\bm{n}_{\bm{c}^{(j\wedge r)}}|\ge1$, one has
\[
\lim_{m\to\infty}\frac{\E[(O_{l,m}^{(n)})^{r}]}{m^{r\sigma}}
=0.
\]
In a similar fashion note that, as $m\to\infty$, the following
asymptotic equivalence holds true:
\begin{eqnarray*}
\E\bigl[\bigl(N_{l,m}^{(n)}\bigr)^{r}\bigr]
&\sim&\Gamma(\theta+n) m^{1-\theta-n}\\
&&{}\times\sum_{y=1}^r
S(r,y)\frac{\sigma^y[(1-\sigma)_{l-1}]^y}{(l!)^y}
\frac{(j+{\theta}/{\sigma})_y}{\Gamma(\theta
+n+y\sigma)}
  m^{\theta+n+y\sigma-1},
\end{eqnarray*}
which, in turn, yields
\[
\lim_{m\rightarrow+\infty}\frac{\E
[(N_{l,m}^{(n)})^{r}]}{m^{r\sigma}}
=\biggl(\frac{\sigma(1-\sigma)_{l-1 }}{l!}\biggr)^{r}
\frac{\Gamma(\theta+n)(j+{\theta}/{\sigma})_{r
}}{\Gamma(\theta+n+r\sigma)}.
\]
Finally, still as $m\to\infty$,
\begin{eqnarray*}
\mathcal{B}^{i}(\sigma,n,j,\bm{n},m)
&\sim&\frac{\Gamma(\theta+n)}{m^{\theta+n-1}}
\sum_{x=1}^{j\wedge i}
x!S(i,x) \sum_{\bm{c}^{(x)}\in\mathcal{C}_{j,x}}
\prod_{t=1}^{x}\frac{(n_{c_{t}}-\sigma
)_{l-n_{c_{t}}}}{(l-n_{c_t})!}\\
&&{}\times
\sum_{y=1}^{r-i} S(r-i,y)\frac{\sigma^{y}[(1-\sigma
)_{l-1}]^y}{(l!)^y}\\
&&{}\times\frac{(j+{\theta}/{\sigma})_{y}}{
\Gamma(\theta+n-|\bm{n}_{\bm{c}^{(x)}}|+x\sigma)}
m^{\theta+n-1+x\sigma-|\bm{n}_{\bm{c}^{(x)}}|}
\end{eqnarray*}
and, since $|\bm{n}_{\bm{c}^{(x)}}|\ge1$ for any $x=1,\ldots,
(j\wedge
i)$, one
has
\[
\lim_{m\to\infty}\frac{1}{m^{r\sigma}} \mathcal{B}^{i}(\sigma
,n,j,\bm
{n},m) =0
\]
for any $i=1,\ldots,r-1$. These limiting relations plainly lead to
conclude that
\begin{eqnarray*}
\lim_{m\rightarrow+\infty}\E\bigl[m^{-r\sigma}
\bigl(M_{l,m}^{(n)}
\bigr)^{r} \bigr]
&=&\biggl(\frac{\sigma(1-\sigma)_{l-1 }}{l!}\biggr)^{r} \frac
{\Gamma
(\theta+n) (j+{\theta}/{\sigma})_{r}}{\Gamma
(\theta
+n+r\sigma)}\\
&=&\biggl(\frac{\sigma(1-\sigma)_{l-1 }}{l!}\biggr)^{r}\E[Z_{n,j}^{r}].
\end{eqnarray*}
According to \cite{Fav09}, Proposition 2, the distribution of the
random variable $Z_{n,j}$ is uniquely
characterized by the moment sequence
$(\E[(Z_{n,j})^{r}])_{r\geq1} $. Similar arguments lead to determine the
limiting distribution of the random variable $N_{l,m}^{(n)}/m_{\sigma
}$, as
$m\rightarrow+\infty$.

\subsection{Proofs for the Gnedin model}\label{sec5.7}

\subsubsection{\texorpdfstring{Proof of Propositions~\protect\ref{prppriorgnedin} and \protect\ref{prpgnedin}}
{Proof of Propositions 8 and 9}}\label{sec5.7.1}
The proof of \eqref{eqpriormomgnedin} follows from
\eqref{eqmomentmarginalprior} and \eqref{eqgnedinweight}, after
noting that $\Ccr(n,k;-1)=(-1)^kn!(n-1)!/[k!(k-1)!(n-k)!]$. As for the
determination of the distributions of $O_{l,m}^{(n)}$ and
$N_{l,m}^{(n)}$, one uses the fact that $\Ccr(n,k;-1,\gamma
)=(-1)^k{{n-\gamma-1\choose n-k}} n!/k!$ along with the results stated
in Theorems~\ref{posteriormarginalold} and \ref
{posteriormarginalnew}.


%

\printaddresses


\begin{thebibliography}{27}

\bibitem{Arr92}
\begin{barticle}[mr]
\bauthor{\bsnm{Arratia},~\bfnm{Richard}\binits{R.}},
  \bauthor{\bsnm{Barbour},~\bfnm{A.~D.}\binits{A.~D.}} \AND
  \bauthor{\bsnm{Tavar{\'e}},~\bfnm{Simon}\binits{S.}}
(\byear{1992}).
\btitle{Poisson process approximations for the {E}wens sampling formula}.
\bjournal{Ann. Appl. Probab.}
\bvolume{2}
\bpages{519--535}.
\bid{issn={1050-5164}, mr={1177897}}
\bptok{imsref}%
\end{barticle}
\endbibitem

\bibitem{Arr03}
\begin{bbook}[mr]
\bauthor{\bsnm{Arratia},~\bfnm{Richard}\binits{R.}},
  \bauthor{\bsnm{Barbour},~\bfnm{A.~D.}\binits{A.~D.}} \AND
  \bauthor{\bsnm{Tavar{\'e}},~\bfnm{Simon}\binits{S.}}
(\byear{2003}).
\btitle{Logarithmic Combinatorial Structures: A Probabilistic Approach}.
\bpublisher{European Mathematical Society}, \baddress{Z\"urich}.
\bid{doi={10.4171/000}, mr={2032426}}
\bptok{imsref}%
\end{bbook}
\endbibitem

\bibitem{Bar92}
\begin{barticle}[mr]
\bauthor{\bsnm{Barbour},~\bfnm{A.~D.}\binits{A.~D.}}
(\byear{1992}).
\btitle{Refined approximations for the {E}wens sampling formula}.
\bjournal{Random Structures Algorithms}
\bvolume{3}
\bpages{267--276}.
\bid{doi={10.1002/rsa.3240030306}, issn={1042-9832}, mr={1164840}}
\bptok{imsref}%
\end{barticle}
\endbibitem

\bibitem{Cha05}
\begin{bbook}[mr]
\bauthor{\bsnm{Charalambides},~\bfnm{Charalambos~A.}\binits{C.~A.}}
(\byear{2005}).
\btitle{Combinatorial Methods in Discrete Distributions}.
\bpublisher{Wiley-Interscience}, \baddress{Hoboken, NJ}.
\bid{doi={10.1002/0471733180}, mr={2131068}}
\bptok{imsref}%
\end{bbook}
\endbibitem

\bibitem{Dur09}
\begin{barticle}[pbm]
\bauthor{\bsnm{Durden},~\bfnm{Chris}\binits{C.}} \AND
  \bauthor{\bsnm{Dong},~\bfnm{Qunfeng}\binits{Q.}}
(\byear{2009}).
\btitle{RICHEST---A web server for richness estimation in biological data}.
\bjournal{Bioinformation}
\bvolume{3}
\bpages{296--298}.
\bid{issn={0973-2063}, pmcid={2655047}, pmid={19293995}}
\bptok{imsref}%
\end{barticle}
\endbibitem

\bibitem{Ewe72}
\begin{barticle}[mr]
\bauthor{\bsnm{Ewens},~\bfnm{W.~J.}\binits{W.~J.}}
(\byear{1972}).
\btitle{The sampling theory of selectively neutral alleles}.
\bjournal{Theoret. Population Biology}
\bvolume{3}
\bpages{87--112}.
\bid{issn={0040-5809}, mr={0325177}}
\bptnote{check related}%
\bptok{imsref}%
\end{barticle}
\endbibitem

\bibitem{Ewe98}
\begin{bincollection}[mr]
\bauthor{\bsnm{Ewens},~\bfnm{W.~J.}\binits{W.~J.}} \AND
  \bauthor{\bsnm{Tavar{\'{e}}},~\bfnm{S.}\binits{S.}}
(\byear{1998}).
\btitle{The Ewens sampling formula, {U}pdate {V}ol. 2}.
In \bbooktitle{Encyclopedia of Statistical Science}
(\beditor{\binits{S.}~\bsnm{Kotz}}, \beditor{\binits{C.~B.}~\bsnm{Read}} \AND \beditor{\binits{D.~L.}~\bsnm{Banks}},
  eds.)
\bpages{230--234}.
\bpublisher{Wiley}, \baddress{New York}.
\bptok{imsref}%
\end{bincollection}
\endbibitem

\bibitem{Fav09}
\begin{barticle}[mr]
\bauthor{\bsnm{Favaro},~\bfnm{Stefano}\binits{S.}},
  \bauthor{\bsnm{Lijoi},~\bfnm{Antonio}\binits{A.}},
  \bauthor{\bsnm{Mena},~\bfnm{Rams{\'e}s~H.}\binits{R.~H.}} \AND
  \bauthor{\bsnm{Pr{\"u}nster},~\bfnm{Igor}\binits{I.}}
(\byear{2009}).
\btitle{Bayesian non-parametric inference for species variety with a
  two-parameter {P}oisson--{D}irichlet process prior}.
\bjournal{J. R. Stat. Soc. Ser. B Stat. Methodol.}
\bvolume{71}
\bpages{993--1008}.
\bid{doi={10.1111/j.1467-9868.2009.00717.x}, issn={1369-7412}, mr={2750254}}
\bptok{imsref}%
\end{barticle}
\endbibitem

\bibitem{Fer73}
\begin{barticle}[mr]
\bauthor{\bsnm{Ferguson},~\bfnm{Thomas~S.}\binits{T.~S.}}
(\byear{1973}).
\btitle{A {B}ayesian analysis of some nonparametric problems}.
\bjournal{Ann. Statist.}
\bvolume{1}
\bpages{209--230}.
\bid{issn={0090-5364}, mr={0350949}}
\bptok{imsref}%
\end{barticle}
\endbibitem

\bibitem{Gne10}
\begin{barticle}[mr]
\bauthor{\bsnm{Gnedin},~\bfnm{Alexander}\binits{A.}}
(\byear{2010}).
\btitle{A species sampling model with finitely many types}.
\bjournal{Electron. Commun. Probab.}
\bvolume{15}
\bpages{79--88}.
\bid{doi={10.1214/ECP.v15-1532}, issn={1083-589X}, mr={2606505}}
\bptok{imsref}%
\end{barticle}
\endbibitem

\bibitem{Gne05}
\begin{barticle}[mr]
\bauthor{\bsnm{Gnedin},~\bfnm{A.}\binits{A.}} \AND
  \bauthor{\bsnm{Pitman},~\bfnm{J.}\binits{J.}}
(\byear{2005}).
\btitle{Exchangeable {G}ibbs partitions and {S}tirling triangles}.
\bjournal{Zap. Nauchn. Sem. S.-Peterburg. Otdel. Mat. Inst. Steklov. (POMI)}
\bvolume{325}
\bpages{83--102}.
\bid{doi={10.1007/s10958-006-0335-z}, issn={0373-2703}, mr={2160320}}
\bptok{imsref}%
\end{barticle}
\endbibitem

\bibitem{griffiths07}
\begin{barticle}[mr]
\bauthor{\bsnm{Griffiths},~\bfnm{Robert~C.}\binits{R.~C.}} \AND
  \bauthor{\bsnm{Span{\`o}},~\bfnm{Dario}\binits{D.}}
(\byear{2007}).
\btitle{Record indices and age-ordered frequencies in exchangeable {G}ibbs
  partitions}.
\bjournal{Electron. J. Probab.}
\bvolume{12}
\bpages{1101--1130}.
\bid{doi={10.1214/EJP.v12-434}, issn={1083-6489}, mr={2336601}}
\bptok{imsref}%
\end{barticle}
\endbibitem

\bibitem{Ho}
\begin{bmisc}[auto:STB|2012/04/30|08:06:40]
\bauthor{\bsnm{Ho},~\bfnm{M.~W.}\binits{M.~W.}},
  \bauthor{\bsnm{James},~\bfnm{L.~F.}\binits{L.~F.}} \AND
  \bauthor{\bsnm{Lau},~\bfnm{J.~W.}\binits{J.~W.}}
(\byear{2007}).
\bhowpublished{Gibbs partitions (EPPF's) derived from a stable subordinator are
  Fox H and Meijer G transforms. MatharXiv preprint. Available at
  arXiv:\arxivurl{0708.0619v2}.}
\bptok{imsref}%
\end{bmisc}
\endbibitem

\bibitem{Jam10}
\begin{barticle}[mr]
\bauthor{\bsnm{James},~\bfnm{Lancelot~F.}\binits{L.~F.}}
(\byear{2010}).
\btitle{Lamperti-type laws}.
\bjournal{Ann. Appl. Probab.}
\bvolume{20}
\bpages{1303--1340}.
\bid{doi={10.1214/09-AAP660}, issn={1050-5164}, mr={2676940}}
\bptok{imsref}%
\end{barticle}
\endbibitem

\bibitem{Kin78}
\begin{barticle}[mr]
\bauthor{\bsnm{Kingman},~\bfnm{J.~F.~C.}\binits{J.~F.~C.}}
(\byear{1978}).
\btitle{The representation of partition structures}.
\bjournal{J. Lond. Math. Soc. (2)}
\bvolume{18}
\bpages{374--380}.
\bid{doi={10.1112/jlms/s2-18.2.374}, issn={0024-6107}, mr={0509954}}
\bptok{imsref}%
\end{barticle}
\endbibitem

\bibitem{Kin82}
\begin{barticle}[mr]
\bauthor{\bsnm{Kingman},~\bfnm{J.~F.~C.}\binits{J.~F.~C.}}
(\byear{1982}).
\btitle{The coalescent}.
\bjournal{Stochastic Process. Appl.}
\bvolume{13}
\bpages{235--248}.
\bid{doi={10.1016/0304-4149(82)90011-4}, issn={0304-4149}, mr={0671034}}
\bptok{imsref}%
\end{barticle}
\endbibitem

\bibitem{lmp07bmc}
\begin{barticle}[mr]
\bauthor{\bsnm{Lijoi},~\bfnm{Antonio}\binits{A.}},
  \bauthor{\bsnm{Mena},~\bfnm{Rams{\'e}s~H.}\binits{R.~H.}} \AND
  \bauthor{\bsnm{Pr{\"u}nster},~\bfnm{Igor}\binits{I.}}
(\byear{2007}).
\btitle{A Bayesian nonparametric method for prediction in EST analysis}.
\bjournal{BMC Bioinformatics}
\bvolume{8}
\bpages{339}.
\bptok{imsref}%
\end{barticle}
\endbibitem

\bibitem{Lij207}
\begin{barticle}[mr]
\bauthor{\bsnm{Lijoi},~\bfnm{Antonio}\binits{A.}},
  \bauthor{\bsnm{Mena},~\bfnm{Rams{\'e}s~H.}\binits{R.~H.}} \AND
  \bauthor{\bsnm{Pr{\"u}nster},~\bfnm{Igor}\binits{I.}}
(\byear{2007}).
\btitle{Bayesian nonparametric estimation of the probability of discovering new
  species}.
\bjournal{Biometrika}
\bvolume{94}
\bpages{769--786}.
\bid{doi={10.1093/biomet/asm061}, issn={0006-3444}, mr={2416792}}
\bptok{imsref}%
\end{barticle}
\endbibitem

\bibitem{lmp07}
\begin{barticle}[mr]
\bauthor{\bsnm{Lijoi},~\bfnm{Antonio}\binits{A.}},
  \bauthor{\bsnm{Mena},~\bfnm{Rams{\'e}s~H.}\binits{R.~H.}} \AND
  \bauthor{\bsnm{Pr{\"u}nster},~\bfnm{Igor}\binits{I.}}
(\byear{2007}).
\btitle{Controlling the reinforcement in {B}ayesian non-parametric mixture
  models}.
\bjournal{J. R. Stat. Soc. Ser. B Stat. Methodol.}
\bvolume{69}
\bpages{715--740}.
\bid{doi={10.1111/j.1467-9868.2007.00609.x}, issn={1369-7412}, mr={2370077}}
\bptok{imsref}%
\end{barticle}
\endbibitem

\bibitem{Lij308}
\begin{barticle}[mr]
\bauthor{\bsnm{Lijoi},~\bfnm{Antonio}\binits{A.}},
  \bauthor{\bsnm{Pr{\"u}nster},~\bfnm{Igor}\binits{I.}} \AND
  \bauthor{\bsnm{Walker},~\bfnm{Stephen~G.}\binits{S.~G.}}
(\byear{2008}).
\btitle{Bayesian nonparametric estimators derived from conditional {G}ibbs
  structures}.
\bjournal{Ann. Appl. Probab.}
\bvolume{18}
\bpages{1519--1547}.
\bid{doi={10.1214/07-AAP495}, issn={1050-5164}, mr={2434179}}
\bptok{imsref}%
\end{barticle}
\endbibitem

\bibitem{Pit95}
\begin{barticle}[mr]
\bauthor{\bsnm{Pitman},~\bfnm{Jim}\binits{J.}}
(\byear{1995}).
\btitle{Exchangeable and partially exchangeable random partitions}.
\bjournal{Probab. Theory Related Fields}
\bvolume{102}
\bpages{145--158}.
\bid{doi={10.1007/BF01213386}, issn={0178-8051}, mr={1337249}}
\bptok{imsref}%
\end{barticle}
\endbibitem

\bibitem{Pit03}
\begin{bincollection}[mr]
\bauthor{\bsnm{Pitman},~\bfnm{Jim}\binits{J.}}
(\byear{2003}).
\btitle{Poisson--{K}ingman partitions}.
In \bbooktitle{Statistics and Science: A~{F}estschrift for {T}erry {S}peed}
(\beditor{\binits{D.~R.} \bsnm{Goldstein}}, ed.).
\bseries{Institute of Mathematical Statistics Lecture Notes---Monograph Series}
\bvolume{40}
\bpages{1--34}.
\bpublisher{IMS}, \baddress{Beachwood, OH}.
\bid{doi={10.1214/lnms/1215091133}, mr={2004330}}
\bptok{imsref}%
\end{bincollection}
\endbibitem

\bibitem{Pit06}
\begin{bbook}[mr]
\bauthor{\bsnm{Pitman},~\bfnm{J.}\binits{J.}}
(\byear{2006}).
\btitle{Combinatorial Stochastic Processes}.
\bseries{Lecture Notes in Math.}
\bvolume{1875}.
\bpublisher{Springer}, \baddress{Berlin}.
\bid{mr={2245368}}
\bptok{imsref}%
\end{bbook}
\endbibitem

\bibitem{quack}
\begin{barticle}[pbm]
\bauthor{\bsnm{Quackenbush},~\bfnm{J.}\binits{J.}},
  \bauthor{\bsnm{Cho},~\bfnm{J.}\binits{J.}},
  \bauthor{\bsnm{Lee},~\bfnm{D.}\binits{D.}},
  \bauthor{\bsnm{Liang},~\bfnm{F.}\binits{F.}},
  \bauthor{\bsnm{Holt},~\bfnm{I.}\binits{I.}},
  \bauthor{\bsnm{Karamycheva},~\bfnm{S.}\binits{S.}},
  \bauthor{\bsnm{Parvizi},~\bfnm{B.}\binits{B.}},
  \bauthor{\bsnm{Pertea},~\bfnm{G.}\binits{G.}},
  \bauthor{\bsnm{Sultana},~\bfnm{R.}\binits{R.}} \AND
  \bauthor{\bsnm{White},~\bfnm{J.}\binits{J.}}
(\byear{2001}).
\btitle{The TIGR gene indices: Analysis of gene transcript sequences in highly
  sampled eukaryotic species}.
\bjournal{Nucleic Acids Res.}
\bvolume{29}
\bpages{159--164}.
\bid{issn={1362-4962}, pmcid={29813}, pmid={11125077}}
\bptnote{check year}%
\bptok{imsref}%
\end{barticle}
\endbibitem

\bibitem{Sch10}
\begin{barticle}[mr]
\bauthor{\bsnm{Schweinsberg},~\bfnm{Jason}\binits{J.}}
(\byear{2010}).
\btitle{The number of small blocks in exchangeable random partitions}.
\bjournal{ALEA Lat. Am. J. Probab. Math. Stat.}
\bvolume{7}
\bpages{217--242}.
\bid{issn={1980-0436}, mr={2672786}}
\bptok{imsref}%
\end{barticle}
\endbibitem

\bibitem{Val09}
\begin{bmisc}[auto:STB|2012/04/30|08:06:40]
\bauthor{\bsnm{Valen},~\bfnm{E.}\binits{E.}}
(\byear{2009}).
\bhowpublished{Deciphering transcriptional regulation---Computational
  approaches. Ph.D. thesis, Bioinformatics Centre, Univ. Copenhagen}.
\bptok{imsref}%
\end{bmisc}
\endbibitem

\end{thebibliography}
\end{document}